\documentclass[a4paper,french,12pt]{amsart}
\usepackage{bourbaki}
\usepackage[all]{xy}
\usepackage[frenchb]{babel}
\usepackage[applemac]{inputenc}
\usepackage{eucal}
\usepackage{mathrsfs}
\usepackage{amssymb}
\usepackage{amsxtra}
\usepackage{enumerate}
\makeatletter
\let\@@seccntformat\@seccntformat
\renewcommand*{\@seccntformat}[1]{%
  \expandafter\ifx\csname @seccntformat@#1\endcsname\relax
    \expandafter\@@seccntformat
  \else
    \expandafter
      \csname @seccntformat@#1\expandafter\endcsname
  \fi
    {#1}%
}
\newcommand*{\@seccntformat@section}[1]{%
  \S\csname the#1\endcsname\quad
}
\makeatother

\makeatletter
\renewcommand{\tocsection}[3]{%
  \indentlabel{\@ifnotempty{#2}{\S~\ignorespaces#1 #2.\quad}}#3}
\makeatother

\addto\captionsfrench{

}

\makeatletter
\@addtoreset{section}{part}% ou \numberwithin{section}{part}
\@addtoreset{proposition}{section}% ou \numberwithin{proposition}{part}
\@addtoreset{theoreme}{part}% ou \numberwithin{theoreme}{part}
\makeatother

\input cyracc.def 
 
\newtheorem{theoreme}{Théorème}
\newtheorem{proposition}{ Proposition}
\newtheorem*{lemme}{ Lemme}
\newtheorem*{corollaire}{Corollaire}

\newtheorem*{definition}{Définition}

\theoremstyle{remark}

\newtheorem{exemple}{\it Exemple}

\def \1{\mathbb {1}}
\def \RM{\mathbb {R}}%        corps des reels
\def \NM{\mathbb{N}}%        entiers naturels
\def \ZM{\mathbb{Z}}%        entiers relatifs
\def \CM{\mathbb{C}}%        nombres complexes

\def \QM{\mathbb{Q}}%        nombres rationnels
\def \PM {\mathbb{P}}

\def \Ker {{\rm Ker\,}}
\def \Im {{\rm Im\,}}

\def \Der {{\rm Der\,}}

\def \Aut {{\rm Aut\,}}

\def \Isom {{\rm Isom\,}}

\def \p {{\rm exp\,}}
\def \Id {{\rm Id\,}}

\def \d{\partial}%derivee partielle
\def\dt{\delta} 
\def\a{\alpha}
\def\b{\beta}
\def\e{\varepsilon}  
\def\g{\gamma}

\def\l{\lambda}

\def\p{\varphi}

\def\G{\Gamma}   
\def\D{\Delta}
\def \s{\sigma}

\def \to{\longrightarrow} 
\def \w{\wedge}
\def \alg{\mathfrak{g}}

\def \< {{\langle }}
\def \> {{\rangle }}
\def \( {\left( }
\def \) {\right) }

\newcommand{\Bt}{{\mathcal B}}

\newcommand{\Et}{{\mathcal E}}
\newcommand{\Ft}{{\mathcal F}}

\newcommand{\Ht}{{\mathcal H}}

\newcommand{\Lt}{{\mathcal L}}
\newcommand{\Mt}{{\mathcal M}}

\newcommand{\Ot}{{\mathcal O}}

\newcommand{\Rt}{{\mathcal R}}

\newcommand{\Vt}{{\mathcal V}}

\renewcommand{\mod}{{\rm  mod\,}}
\title[Espaces vectoriels \'echelonn\'es]{Espaces vectoriels \'echelonn\'es}
\author{  Mauricio  Garay}
\address{Max Planck,\ Institut für Mathematik\\
Vivatsgasse 7\\
53111 Bonn, Allemagne.}
\email{garay@mpim-bonn.mpg.de}
 
\begin{document}
\maketitle
%\begin{flushright}{\em Cet article est dédié à Duco van Straten.}\end{flushright}
%\bigskip
\section*{Introduction}
Dans l'étude des singularités d'applications différentiables, les relations entre les orbites de l'action du groupe des difféomorphismes et celles de son action infinitésimale sont décrites par les théorèmes de déformations verselles, de stabilité par déformation et de détermination finie~ \cite{AVG,Mather_fdet,Mather_stability,Tougeron,Tyurina}. Le but de cette article est d'étendre ces résultats lorsque des petits dénominateurs interveniennent  dans les équations cohomologiques.

Historiquement, de telles recherches sont, en fait, bien antérieures à celles effectuées en théorie des singularités.
Elles ont vu le jour avec la première thèse de Poincaré, pour aboutir dans les années 50  au théorème des tores invariants de Kolmogorov. Ce dernier, en modifiant l'algorithme de Newton a ouvert la voie à ce qu'Arnold appelait la lutte contre les petits dénominateurs~\cite{Arnold_SD1,Arnold_edo,Kolmogorov_KAM}. 

De son côté, Moser a préféré développer une approche basée sur les théorèmes de fonctions implicites~\cite{Moser_KAM}  (voir également \cite{Hamilton_implicit,Nash_imbedding}). Il en a déduit une heuristique sur les actions de groupes en dimension infinie~\cite{Moser_pde}. Le problème posé par cette approche n'avait pas échappé à son auteur~: elle ne donnait pas la relation souhaitée entres l'action du groupe et l'action infinitésimale, mais seulement avec l'action linéarisée en chaque point du groupe. Certes dans un groupe $G$, la différentielle de la multiplication par $g \in G$ envoie l'algèbre de Lie du groupe sur le plan tangent au groupe en $g$ mais, en dimension infinie, cette opération est la source de nombreuses difficultés.

 Pour tenter remédier à ce problème, Moser suggéra que l'on pouvait se contenter d'un inverse approché pour l'action. Dans thèse, Sergeraert
 construisit un cadre formel pour l'heuristique de Moser. Il réussit, dans certains cas, à controller le problème posé par la multiplication dans le groupe~\cite[Théorème 4.2.5 et Corrollaire 4.2.6]{Sergeraert}.  Mais, de nombreux problèmes demeuraient inaccessible par l'approche de Sergeraert. 
  
C'est ce que vit Zehnder qui formalisa la notion d'inverse approché. Mais à l'instar de ses précédécesseurs, il se voyait forcé d'admettre qu'il n'était parvenu à démontrer aucun résultat essentiellement nouveau, tout au plus avait-il amélioré les conditions de différentiabilité~: «We illustrate our method by small divisor theorems, which are well-known, but our proof is in principle simpler and more systematic»~\cite[Introduction]{Zehnder_implicit}. Dans le contexte général des actions de groupes, Zehnder  fût contraint  de reproduire une heuristique semblable à celle de Moser~\cite[Chapter 5]{Zehnder_implicit}.  
 
Dans les années 80, Herman parvint à une nouvelle démonstration  du théorème KAM, grâce aux théorèmes de  fonctions implicites
obtenus par Hamilton, Sergeraert et Zehnder~\cite{Bost_KAM}.  Cette démonstration le conduisit aux premiers exemples de théorèmes de type KAM pour des variétés involutives de dimension arbitraire~\cite{Yoccoz_Herman}. Herman obtenait ainsi des obstructions à l'hypothèse ergodique.
 
C'est probablement, les relations entre les tores invariants et l'hypothèse ergodique qui amenèrent Herman à formuler la conjecture suivante~\cite{Herman_ICM}:\\
  
 {\em  Au voisinage d'un point fixe elliptique dont le linéarisé est diophantien, tout symplectomorphisme réel analytique possède un ensemble de mesure positive de tores invariants.}\\
 
 Cet article  a pour but de poser des fondations qui permettent, entre autres, de résoudre cette conjecture. 
 
Pour cela, je me propose d'initier une théorie générale des actions de groupes en dimension infinie dans le contexte de la géométrie analytique.

 Cet article est organisé de la façon suivante:\\
  Au §1, on définit  la catégorie des espaces vectoriels échelonnés.  Dans les théories différentiables, comme celles de Sergeraert et de Zehnder, l'espace fonctionnel est représenté comme une limite inverse d'espaces de Banach, ce qui le munit d'une structure d'espace de Fréchet.
  Dans la théorie analytique, on peut choisir entre limite inverse et limite directe. Nous ferons le choix de la limite directe bien que celle-ci soit, en général, non métrisable.  En algèbre, un module peut avoir de nombreuses résolutions ou en topologie, une variété peut avoir de nombreuses décompositions cellulaires. Il en va de même pour l'échelonnement d'un espace vectoriel, de nombreuses possibilités s'offrent naturellement et on considère tous ces choix, en tentant d'introduire des notions qui ont de bonnes propriétés fonctorielles.

  Au §2, je rappelle le résultat classique de Mather et Tougeron sur la détermination finie des germes dans un cadre holomorphe, puis le théorème de Siegel-Arnold sur la linéarisation des germes de fonctions holomorphes. Je donne à ces deux théorèmes une formulation analogue afin d'introduire
  le lecteur à la notion de détermination finie d'un point de vue abstrait.  On généralise ensuite la filtration de l'anneau des germes de fonctions par les puissances de l'idéal maximal. Je termine le § par un rappel sur les différents types d'algorithmes itératifs.
  
  Au §3, on démontre la convergence du procédé itératif de Kolmogorov mis en place au § précédent.  Je donne un critère simple de convergence d'un produit infini de morphismes~(Théorème \ref{T::produits}). Bien qu'élémentaire, ce résultat est la clé de voûte de l'édifice~: une fois que l'on dispose d'un critère simple de convergence, il n'y aucune difficulté à démontrer la convergence de procédés itératifs. On démontre ainsi le théorème de $M$-détermination, généralisation des résultats de Mather-Tougeron et de Arnold-Siegel.
  
Au §4, on généralise le théorème de $M$-détermination au cas où les objets sont déterminés modulo une transversale.  

  Au §5, j'introduis les déformations d'espaces vectoriels échelonnés.  Ces considérations sont nécessaires pour l'étude des systèmes dynamiques comme le montreront les exemples du §7.
 
  Au §6, je montre que les structures échelonnées existent naturellement en géométrie analytique. Il est probable qu'elles pourraient éviter, dans certains cas, d'avoir recours aux voisinages privilégiés et, par exemple, de simplifier les arguments de la thèse Douady pour la construction de l'espace des modules des sous-variétés complexes~\cite{Douady_these}.
  
  Au §7, je montre comment les conditions diophantiennes habituelles permettent d'appliquer les théorèmes de $M$-détermination et de transversalité. Je termine le §, par trois applications à la théorie KAM. D'un point de vue abstrait, cette théorie apparaît comme une variante symplectique de l'étude de fonctions sur les variétés (voir e.g. \cite{AVGL2,Goryunov}). Je donne un théorème des tores invariants sans condition de
  non-dégénérescence.  Ce résultat est à rapprocher de celui que j'avais annoncé  en 2009, concernant l'inutilité  des conditions de non-dégénérescence pour garantir l'existence d'un ensemble de mesure positive de tores invariants. Il est également proche du travail (non-publié) de Eliasson-Fayad-Krikorian relié à la conjecture de Herman. 
  
Ensuite,  je montre qu'en un point critique de Morse, le flot hamiltonien est conjugué à un flot linéaire sur une variété lagrangienne complexe singulière pourvu que la fréquence soit diophantienne. Dans le cas réel elliptique, l'origine est le seul point réel de cette variété et ce résultat est donc vide. C'est ce qui explique, peut-être, que les géomètres l'aient ignoré jusqu'à aujourd'hui.  Il ne s'agit pourtant  pas une simple curiosité, puisqu'en quantifiant la conjugaison, on montre que la série perturbative du spectre de l'état fondamental est de classe Gevrey 2 dans la constante de Planck.  
  
   \noindent {\bf Remerciements.} { Pendant plusieurs années, les résultats de cet article sont demeurés à l'état de conjectures. Je remercie J. Féjoz.  Sans son aide,  ces conjectures n'auraient très probablement jamais accédé au rang de théorèmes. Je remercie également H. Eliasson qui m'a aidé à faire mes premiers pas en théorie KAM, ainsi que R. Krikorian et B. Fayad pour leurs  explications sur la conjecture de Herman.}
\tableofcontents
%%%%%%%%%%%%%%%%%%%%%%%%%%%%
 \section{La catégorie des espaces vectoriels échelonnés}
  %%%%%%%%%%%%%%%%%%%%%%%%%%%%%%
  %%%%%%%%%%%%
 %%%%%%%%%%%%%%%%%%%%%%%%%%%%%
\subsection{Définition}
\label{SS::definition}
Dans tout ce qui va suivre, on peut sans difficulté considérer le cas d'espaces vectoriels sur $\CM$ plutôt que sur $\RM$.

Une {\em $S$-échelle de Banach} est une famille décroissante d'espaces de Banach $(E_s)$, $s \in ]0,S[$, telle que 
les inclusions 
$$E_{s+\s} \subset E_s,\ s \in ]0,S[,\ \s \in ]0,S-s[$$ 
soient de norme au plus~$ 1$.

Soit $E$ un espace vectoriel topologique. Un {\em $S$-échelonnement} de $E$ est une échelle  $(E_s)$
de sous-espaces de Banach de $E$ telle que
\begin{enumerate}[{\rm i)}] 
\item  $E=\bigcup_{s \in ]0,S[} E_s$ ;
\item la topologie  limite directe de la topologie des espaces de Banach $E_s$ coïncide avec celle de $E$.
\end{enumerate}
(Rappelons que si $f_s:X_s \to X$ un famille d'application d'espaces topologiques $(X_s)$ dans un ensemble $X$. On appelle topologie {\em  limite directe} sur $X$, la topologie  la plus fine sur $X$ qui rend les applications $f_s$ continues~:
 $$U \subset X {\rm \ est\ ouvert\ }\ \iff  f^{-1}_s(U) {\rm \  est\ ouvert\  dans}\ X_s,\ {\rm pour\ tout}\ s.)$$
 
 L'intervalle $]0,S[$ s'appelle {\em l'intervalle d'échelonnement}.  Si $F$ est un sous-espace vectoriel fermé d'un espace vectoriel échelonné $E$ alors $E/F$ est échelonné par les espaces de Banach
$E_s/(E \cap F)_s$.
  
La notion d'échelonnement vise à transférer les propriétés de $E$ aux espaces de Banach $E_s$, mais l'objet que l'on étudie reste $E$ et non pas l'échelle de Banach. Lorsque le paramètre $S$ ne joue pas de rôle particulier, nous parlerons simplement d'échelle de Banach ou d'espace vectoriel échelonné. 

L'utilisation d'échelles de Banach en analyse remonte aux fondements de l'analyse fonctionnelle. On la trouve par exemple dans la démonstration du théorème de Cauchy-Kovalevskaïa donnée en 1942 par Nagumo~\cite{Nagumo} (voir également \cite{Ovsyannikov}).
 Elle est également à la base de la démonstration proposée par Kolmogorov du théorème des tores invariants~\cite{Kolmogorov_KAM}. 
 
 Cependant ces auteurs ne considèrent qu'une échelle fixe, l'idée de considérer toutes les échelles possibles d'un sous-espace vectoriel topologique est déjà présente dans la thèse de Grothendieck~\cite{Grothendieck_PTT}. En revanche, Grothendieck n'utilise pas le choix d'un paramétrage de l'échelle comme une donnée supplémentaire.
 
   %%%%%%%%%%%%%%%%
 \subsection{L'espace des polynômes}
 
Voici un exemple  simple d'espace vectoriel échelonné. Nous n'utiliserons pas cet exemple dans les applications concrètes, mais il est utile pour comprendre les notions générales, au même titre que $\RM^2$ est un espace de Hilbert qui a son importance.

  Considérons l'espace vectoriel $\RM[X]$  des polynômes en une indéterminée sur $\RM$. Munissons cet espace de la topologie pour laquelle les ouverts sont les ensembles dont l'intersection avec tout sous-espace vectoriel de dimension finie est ouvert
$$U{\rm\ ouvert\ de\ } \RM[X] \iff U \cap F {\rm\ ouvert\ de\ } F,\ \dim F<+\infty .$$
L'espace $\RM[X]$ est l'union des espaces $\RM_k[X]$ des polynômes de degré au plus $k$.
Identifions l'espace $\RM_k[X]$ à $\RM^{k+1}$ par l'application
$$a_k X^k+\cdots+a_0 \mapsto (a_k,\dots,a_0). $$
La norme euclidienne de $\RM^{k+1}$ induit sur $\RM_k[X]$ une structure d'espace vectoriel normé. Pour faire de $\RM[X]$ un espace vectoriel échelonné, il faut fixer la relation entre le paramètre $s$ de l'échelonnement et le degré des polynômes.
 
 Faisons par exemple le choix suivant~: prenons pour $E_s$, l'espace des polynômes dont
le degré est au plus la partie entière de $1/s$. On obtient ainsi une structure d'espace vectoriel échelonné sur $\RM[X]$.

 Plus généralement, en partant d'un espace vectoriel gradué
$$ V=\oplus_{n \geq 0} V_n,\ \dim V_n < +\infty,$$
on définit alors une structure échelonnée sur  $V$ en posant $E_s=\oplus_{n \leq 1/s} V_n$.

% %%%%%%%%%%%%%%%%%%%%%%%%%%%%%%%%%%%%%%%%%
\subsection{Morphismes d'un espace vectoriel échelonné}
\label{SS::morphismes}
Soit $E,F$ deux espaces vectoriels respectivement $S$-échelonné et $S'$ échelonné.

Nous dirons d'une application linéaire que c'est  un {\em morphisme} entre des espaces vectoriels échelonnés $E,F$, si pour tout $s' \in ]0,S[$, 
il existe $s \in ]0,S[$ tel que l'espace de Banach $E_{s'}$ est envoyé continûment dans $F_s$. Nous avons ainsi définit la  {\em catégorie des espaces vectoriels échelonnés}.

Nous désignerons par $\Lt(E,F)$ l'espace vectoriel des morphismes de $E$ dans $F$ et lorsque $E=F$, nous utiliserons la notation $\Lt(E)$ au lieu de $\Lt(E,E)$. Il n'y pas de raison, a priori, pour que $\Lt(E,F)$ coïncide avec l'espace des applications linéaires continues de $E$ dans $F$, mais dans les exemples concrets que nous allons traiter ce sera toujours le cas.
 
   Si $\| \cdot \|$ désigne la norme d'opérateur sur l'espace de Banach $\Lt(E_{s'},F_s)$, nous noterons $\| u \|$ la norme de l'opérateur défini par restriction de $u$ à $E_{s'}$.

Le noyau d'un  morphisme $ u:E \to F$ entre espaces vectoriels $S$-échelonnés est un espace vectoriel $S$-échelonné par~:
$$ (\Ker u)_s=E_s \cap \Ker u,\ s \leq S.$$
 
Venons-en à la notion de convergence d'une suite de morphismes. La norme d'opérateur induit sur  les espaces vectoriels $\Lt(E_{s'},F_s)$, une structure d'espace de Banach. 
\begin{definition}Une suite de morphismes $(u_n)$ de $\Lt(E,F)$ converge vers un morphisme $u \in \Lt(E,F)$ si pour tout $s' \in ]0,S[$, il existe $s \in ]0,S[$ tel que la restriction de $(u_n)$ définisse une suite de $\Lt(E_{s'},F_s)$ qui converge vers la restriction de $u$.
\end{definition}
Un sous-ensemble $X$ de $\Lt(E,F)$ sera dit {\em fermé} si toute suite convergente de points de $X$ à sa limite dans $X$. (L'utilisation du mot «fermé» est légèrement abusive, car il ne s'agit pas a priori du complémentaire d'un ouvert.)
 %%%%%%%%%%%%%%%
 \subsection{  $\tau$-morphismes, morphismes bornés}

Conservons les notations du chapitre précédent. 
\begin{definition} Un morphisme $ u \in \Lt(E,F)$  est appelé  un $\tau$-morphisme
si pour tout $s'  \in ]0,\tau]$ et pour tout $s \in ]0,s'[$, on a l'inclusion $ u(E_{s'}) \subset F_s$ et $ u$ induit par restriction une application linéaire continue
$$u_{s',s}~:~E_{s'}~\to~F_s.$$
\end{definition}
On a alors des diagrammes commutatifs
$$\xymatrix{ & \ F_s \ar@{^{(}->}[d]\\
 E_{s'} \ar[r]^-{u_{\mid E_{s'}}} \ar[ru]^{u_{s',s}}  & F }
$$
pour tout $s'  \in ]0,\tau]$ et pour tout $s \in ]0,s'[$, la flèche verticale étant donnée par l'inclusion $F_s \subset F$.
\begin{exemple} Reprenons l'exemple du \ref{SS::definition}. Les morphismes sont les applications linéaires qui envoient les sous-espaces de dimension finie sur des sous-espaces de dimension finie. Toute application linéaire est donc un morphisme. On vérifie facilement que cette propriété caractérise les applications linéaires continues de
$\RM[X]$ donc, à l'instar de la dimension finie, toute application linéaire de l'espace des polynômes dans lui-même est continue. Les $\tau$-morphismes sont les morphismes qui envoie $\RM_k[X]$ dans lui-même
pour tout $\displaystyle{k \geq \left[\frac{1}{\tau}\right]}$. La multiplication par $X$ est un exemple de morphisme qui n'est pas un $\tau$-morphisme.
\end{exemple}
\begin{definition}
\label{D::borne}
  Un $\tau$-morphisme $ u:E \to F$ d'espaces vectoriels $S$-échelonnés est dit 
    $k$-borné, $k \geq 0 $ s'il existe un réel $C>0$ tel que~:
  $$| u(x) |_s \leq C \sigma^{-k} | x |_{s+\sigma},\ {\rm pour\ tous\ }\ s \in ]0,\tau[,\ \s \in ]0,\tau-s],\ x\in E_{s+\s} . $$
\end{definition}
Un morphisme est dit {\em $k$-borné} (resp. {\em borné}) s'il existe $\tau$ (resp. $\tau$ et $k$) pour lequel (resp. lesquels) c'est un $\tau$-morphisme $k$-borné. Lorsque $E=E_s$ et $F=F_s$ sont des espaces de Banach, on retrouve la définition habituelle de morphismes bornés.
(Nous n'utiliserons pas la notion plus générale d'application linéaire bornée d'un espace localement convexe, notre terminologie ne devrait donc pas porter à confusion.)  
\begin{exemple}L'application qui à un polynôme associe sa dérivée est un morphisme qui n'est pas borné. Tout morphisme borné de $\RM[X]$ est automatiquement $0$-borné. Il y a donc très peu de morphismes bornés dans $\RM[X]$.
\end{exemple}

L'espace vectoriel des $\tau$-morphismes (resp. des morphismes)  $k$-bornés
entre $E$ et $F$ sera noté $\Bt^k_\tau(E,F)$ (resp. $\Bt^k(E,F)$). On note $N_\tau^k(u)$ la plus petite constante $C$ vérifiant l'inégalité de la définition \ref{D::borne}.

La propriété pour un endomorphisme {\em surjectif} d'être borné passe au quotient~: tout endomorphisme $k$-borné surjectif $ u:E \to E$ définit un endomorphisme $k$-borné sur l'espace quotient $E/F$, pour tout sous-espace vectoriel fermé $F \subset E$.
%%%%%%%%%%%%%%%%%%%%%%%%%%%%%%%%%%%
\subsection{L'échelle de Banach $(\Bt^k_\tau(E,F))$.}
\begin{proposition}[Féjoz] Si $E,F$ sont des espaces vectoriels $S$-échelonnés alors les espaces vectoriels normés
$(\Bt^k_\tau(E,F),N_\tau^k),\ \tau \in ]0,S[$, forment une $S$-échelle de Banach.  
\end{proposition}
 \begin{proof}
  La seule difficulté consiste à montrer que l'espace vectoriel $\Bt^k_\tau(E,F)$ est complet pour la norme $N_\tau^k$, pour tout $\tau \in]0,S[$.
  
 Je dis que toute suite de Cauchy  $(u_n) \subset \Bt^k_\tau(E,F)$  converge vers un morphisme $u \in \Lt(E,F)$ au sens de \ref{SS::morphismes}.
   Soit donc $s' \in ]0, \tau[$ et $s \in ]0,s'[$. Comme les $(u_n)$ sont des $\tau$-morphismes, ils induisent, par
 restriction, des applications linéaires continues 
 $$v_n:E_{s'} \to F_s.$$  
 Par définition de la norme $N_\tau^k$, la suite $(v_n)$ est de Cauchy dans l'espace de Banach $\Lt(E_{s'},F_s)$ donc convergente. Ceci démontre l'affirmation.  
   
 Montrons à présent que si une suite de Cauchy  $(u_n) \subset \Bt^k_\tau(E,F)$  converge vers un morphisme $u \in \Lt(E,F)$ alors
 $u$ est dans $ \Bt^k_\tau(E,F)$. L'inégalité
  $$| N_\tau^k( u_n) - N_\tau^k( u_m) | \leq  N_\tau^k( u_n- u_m) $$
   montre que  la suite $(N_\tau^k(u_n))$ est de Cauchy dans $\RM$ donc majorée par un constante $C>0$.
   On a alors les inégalités~:
   $$| u(x) |_s \leq | u(x)-u_n(x) |_s+C\s^{-k}| x |_{s+\s}, \ {\rm\ pour\ tout\ } n, $$ 
  pour tout $x \in E_s$, pour tout $s \in ]0,\tau]$ et pour tout $\s \in ]0,\tau-s]$. Par conséquent, le $\tau$-morphisme limite $ u$ est $k$-borné de norme au plus égale à $C$.  
\end{proof}
La suite $(\Bt^k_\tau(E,F),N_\tau^k),\ \tau \in ]0,S[$ munit l'espace vectoriel $\Bt^k(E,F)$ d'une structure d'espace vectoriel échelonné. 
%%%%%%%%%%%%%%%%%%%%%%%%%%%%%
\subsection{Produits de morphismes bornés}
\label{SS::produits}
Conservons les notations du n° précédent. Si $u,v$ sont des morphismes,  respectivement $k$ et $k'$ borné, alors leur composition $u  v$ est
  $(k+k')$-borné et on a l'inégalité
  $$N_\tau^{k+k'}( u v) \leq 2^{k+k'} N_\tau^k(u)N_\tau^{k'}(v) .$$
  En effet~:
 $$| (u v)(x) |_s \leq N_\tau^k(u) \frac{2^k}{\sigma^{k}} |v (x)|_{s+\sigma/2} \leq   N_\tau^k(u)N_\tau^{k'}(v) \frac{2^{k+k'}}{\sigma^{k+k'}} |x|_{s+\sigma} $$
 pour  tout $x\in E_{s+\s}$. Plus généralement, on a la
\begin{proposition}
\label{P::morphismes}
Le produit de $n$ morphismes $k_i$ bornés $u_i,\ i=1,\dots,n$, est un morphisme $k $-borné avec $k:=\sum_{i=1}^n k_i$ et plus précisément
  $$N_\tau^k( u_1  \cdots  u_n) \leq n^k \prod_{i=1}^n N_\tau^{k_i}(u_i),\ .$$
  De plus, si $u_1=\dots=u_n=u$ est d'ordre $1$ alors
   $$\frac{N_\tau^n( u^n )}{n!} \leq 3^n  N_\tau^1(u)^n. $$
 \end{proposition} 
  La première partie de la proposition s'obtient de façon analogue au cas $n=2$, en découpant l'intervalle $[s,s+\s]$ en $n$ parties égales.
   En prenant tous les $u_i$ égaux et d'ordre $1$, on obtient alors
  $$N_\tau^n( u^n ) \leq n^n  N_\tau^1(u)^n.$$ 
  Il nous suffit donc de démontrer l'inégalité
  $$3^n n! \geq n^n  $$
Pour cela, on utilise l'expression intégrale suivante de  la fonction $\G$~:
  $$n!=\G(n+1) =n^{n+1} \int_{0}^{+\infty} e^{n (\log t-t)} dt .$$
  En utilisant l'estimation
  $$ \log t-t \geq -\frac{1}{2}(t-1)^2-1$$
  il vient~:
  $$\G(n+1) \geq  n^{n+1}e^{-n} \int_{0}^{+\infty} e^{-\frac{1}{2}(t-1)^2} dt \geq  n^{n+1}e^{-n} \geq n^n 3^{-n}. $$
 (Cette estimation, comme toutes celles qui suivront, seront toujours loin d'être optimales.)
 Ceci démontre la proposition.
 
 \begin{corollaire}
 Soit $u$ un $\tau$-morphisme $1$-borné. Si l'inégalité $3N_s^1(u) ~<~s $ est satisfaite pour tout $s \leq \tau$ alors
 la série 
 $$e^u:=\sum_{j \geq 0}\frac{u^j}{j!} $$
 converge vers un morphisme de $E$, et plus précisément
$$| e^u x |_{\l s} \leq \sum_{j \geq 0} \frac{(3N^1_s(u))^j }{(1-\l)^j s^j} | x |_s=\frac{1}{1-\frac{3N^1_s(u) }{(1-\l) s} }   |x |_s$$
pour tous $\l \in ]0, 1-\frac{3N^1_s(u)}s[$, $s \in ]0, \tau]$ et $x \in E_s$.  
\end{corollaire}
\begin{proof}
On a 
$$| e^u x |_{\l s} \leq \sum_{j \geq 0} \frac{1}{j!}| u^j x |_{\l s}. $$
D'après la proposition précédente, le morphisme $u^j$ est $j$-borné et on a l'inégalité~:
$$ \frac{1}{j!} | u^j x |_{\l s} \leq \frac{(3N^1_s(u))^j }{(1-\l)^j s^j} | x |_s, $$
ce qui démontre le corollaire.
\end{proof}
Un morphisme est dit exponentiable si son exponentielle définit une série convergente au sens de \ref{SS::morphismes}.
Introduisons {\em la condition (E)} pour un $\tau$-morphisme $u$~:\\

\noindent (E) Pour tout $s \leq \tau$, le $\tau$-morphisme $u$ vérifie l'inégalité $3N_s^1(u) ~<~s $.\\

La corollaire précédent montre que si $u$ satisfait la condition $(E)$ alors il est exponentiable. Finalement, remarquons que deux morphismes 1-bornés $u,v \in \Bt^1(E)$ exponentiables qui commutent vérifient l'égalité
 $$e^{u+v}=e^{u}e^{v}.$$ En effet, si $u$ et $v$ commutent alors les suites
 $$A_n:=\sum_{j= 0}^n \frac{u^j}{j!},\ B_n=\sum_{j= 0}^n \frac{v^j}{j!}. $$
 vérifient
 $$A_nB_n=\sum_{j= 0}^n \frac{(u+v)^j}{j!} $$
et si des suites de morphismes $(A_n), (B_n)$ convergent respectivement vers $A,B$ alors $(A_nB_n)$ converge vers $AB$. Le cas particulier $v=-u$ montre que l'exponentielle d'un morphisme $1$-borné est inversible. 
 %%%%%%%%%%%%%%%%%%%%%%%%%%%%%
  %%%%%%%%%%%%%%%%%%%%%%%%%%%%%%%%%%
 \subsection{\'Echelonnement de $\Ot_{\CM,0}$}
%Considérons l'espace vectoriel $\Ot_{\CM,0}$ des germes de fonctions holomorphes à l'origine. Munissons cet espace de la topologie de la convergence uniforme sur les compacts de $\CM$.
Considérons l'espace vectoriel  $\Ot_{\CM,0}$ des germes en l'origine de fonctions holomorphes d'une variable. Nous désignerons également cet espace par $\CM\{ z \} $ lorsque nous voudrons préciser le choix d'une coordonnée. Pour chaque compact $K$ contenant l'origine, on note
$B(K)$ l'espace des fonctions continues sur $K$ et holomorphes dans l'intérieur de $K$. C'est un espace de Banach pour la norme
$$B(K) \to \RM,\ f \mapsto \sup_{z \in K} | f(z)|. $$
(La complétude est une conséquence immédiate de la formule de Cauchy.)

Ordonnons l'ensemble $\Vt$ des voisinages compacts de l'origine par l'inclusion.
La limite directe des $B(K), K \in \Vt $ s'identifie alors avec l'espace vectoriel  $\Ot_{\CM,0}$~:
$$\Ot_{\CM,0}=\underrightarrow{\lim}\,B(K),\ K \in \Vt.$$
L'espace vectoriel $\Ot_{\CM,0}$ se voit ainsi muni d'une structure d'espace vectoriel topologique.

Soit $(K_s), K_s \in \Vt,\ s \in ]0,S[$, un système fondamental de voisinages compacts de l'origine. La suite $(B(K_s))$ définit un échelonnement de   l'espace vectoriel $\Ot_{\CM,0}$. C'est l'exemple le plus classique  d'échelonnement de $\Ot_{\CM,0}$~\cite{Douady_these,Grauert_Remmert,Nagumo}.  

Un résultat dû à Grothendieck montre que tout borné de $\Ot_{\CM,0}$ est contenu dans l'un des $B(K_s)$
(voir \cite{Grothendieck_EVT}, Chapitre 3). De ce résultat, on déduit aisément que l'espace vectoriel des morphismes de $\Ot_{\CM,0}$ dans lui-même coïncide avec celui des applications linéaires continues. De même, la notion de convergence sur $\Lt(\Ot_{\CM,0})$ établie au \ref{SS::morphismes} coïncide avec celle de la topologie forte, les ensembles fermés au sens de \ref{SS::morphismes}
sont les ensembles fermés pour cette topologie.

Les suites $(B(K_{s^2}))$, $(B(K_{3s}))$ donnent d'autres exemples d'échelonnement de $\Ot_{\CM,0}$, avec les mêmes espaces de Banach,
mais indexés de façon différente.

Prenons à présent pour compact $K_s$, le disque fermé $D_s$ de rayon $s$ centré en l'origine et posons
$$E_s=B(D_s),\ | f |_s=\sup_{z \in D_s} | f(z) |.$$
L'application linéaire continue
$$\Ot_{\CM,0} \to \Ot_{\CM,0},\ [z \mapsto f(z)] \mapsto  [z \mapsto f(2z)] $$
 donne un exemple simple d'un morphisme qui n'est pas un $\tau$-morphisme.
% Un théorème classique de Montel permet d'affirmer que c'est une structure échelonnée localement compacte \cite{Montel}. On le voit d'ailleurs facilement, car l'application de restriction $E_{s+\s} \to E_s$
%est limite des opérateurs de rang fini 
%$$E_{s+\s} \to E_s,\ \sum_{n \geq 0} a_n z^n \mapsto \sum_{n \leq N} a_n z^n.$$
 
Je dis que pour l'échelonnement $(E_s)$,  tout opérateur différentiel d'ordre $k$ sur $\Ot_{\CM,0}$ définit un morphisme $k$-borné. D'après la proposition \ref{P::morphismes}, il suffit de montrer que l'application linéaire~:
  $$\d_z:\Ot_{\CM,0} \to \Ot_{\CM,0},\ f \mapsto f' $$
  est $1$-borné. Soit donc $f \in E_{s+\s}$, $z \in D_s$ et notons $\g \subset D_{s+\s}$ le cercle centré en $z$ de rayon $\s$ orienté positivement, après une intégration par parties, la formule de Cauchy donne~: 
   $$f'(z)=\frac{1}{2i\pi}\int_{\g} \frac{f(\zeta)}{(z-\zeta)^2}\, d\zeta ,$$
 d'où l'estimation
  $$| f'(z) | \leq \sigma^{-1} | f |_{s+\sigma},\ {\rm pour\ tout\ } z\in D_s. $$
  Ce qui démontre l'affirmation. 
   
  Remarquons au passage que la dérivée $k$-ième est donnée par la formule
   $$f^{(k)}(z)=\frac{k!}{2i\pi}\int_{\g} \frac{f(\zeta)}{(z-\zeta)^{k+1}}\, d\zeta. $$
  On obtient ainsi les {\em inégalités de Cauchy}  
 $$| f^{(k)}(z) | \leq k!\sigma^{-k} | f |_{s+\sigma},$$
 qui montrent que l'estimation de la proposition \ref{P::morphismes} n'est pas optimale.
 
 Il existe, bien entendu, des opérateurs $k$-bornés qui ne sont pas différentiels, par exemple celui qui envoie $z^k$ sur $k z^{k-1}$ pour $k>0$ et $1$ sur lui même. 
 
 Nous voyons sur cet exemple qu'échelonner l'espace vectoriel topologique $\Ot_{\CM,0}$
 consiste à mettre en relation la taille des voisinages dans $\CM$ (ici le rayon des disques) avec l'indexation de la suite d'espaces de Banach $(E_s)$.
 Cette mise en relation de grandeurs de nature différente constitue le principe  de l'échelonnement. Les échelonnements ont une tendance naturelle à se propager~: un fois mis en rapport l'indice du disque $D_s$ avec son rayon, les espaces de Banach se trouvent également indexés en fonction de ce rayon. Puis l'espace des applications bornées se trouve lui-même échelonné.
  
 Venons en à l'exponentielle.

L'exponentielle de $\d_z$ diverge. Cependant, si $f$ est holomorphe dans un disque de rayon supérieur à $1$, la suite de fonctions
 $(\sum_{j=0}^n \frac{\d_z^j}{j!} f)_{n \in \NM} $ a pour limite le germe  $ [z \mapsto f(z+1)] \in \Ot_{\CM,0}.$

L'exponentielle de $\d_z$ est le flot au temps $t=1$ du champ de vecteurs $\d_z$.  C'est un opérateur non borné, dont le domaine de définition est le sous-espace vectoriel des germes qui sont holomorphes dans un disque de rayon supérieur à $1$. 
 
 L'exponentielle de $\l z \d_z$ converge et donne l'automorphisme
 d'algèbre
 $$z \mapsto (\sum_{n \geq 0} \frac{(\l z\d_z)^n}{n!} )z=e^\l z.  $$
 Plus généralement, toute dérivation de la forme 
$$zh(z)\d_z,\ h\in \Ot_{\CM,0} $$
est exponentiable. En résumé, l'algèbre de Lie $\alg$ des dérivations de $\Ot_{\CM,0}$ est filtrée
$$\alg \supset \Mt^1_{\alg} \supset \Mt^2_{\alg} \supset \cdots  $$
avec
$$ \Mt^k_{\alg}=\{ z^kh(z)\d_z,\ h \in \Ot_{\CM,0} \}$$
et tout élément de $\Mt^1_{\alg}$ est exponentiable.

Notons $\Mt_{\CM,0}$, l'idéal maximal de l'anneau local $\Ot_{\CM,0}$~:
$$\Mt_{\CM,0}=\{ f \in \Ot_{\CM,0}: f(0)=0 \}.$$
 Lorsque $v \in \Mt^2_{\alg}$ et $f \in\Mt_{\CM,0}^k$,on a~:
 $$ e^v f=f+v \cdot f\ (\mod \Mt_{\CM,0}^{k+1}). $$
 Ce qui donne un sens précis au fait que l'action infinitésimale de $e^v$ est donnée par la dérivation le long de $v$. En général, ce n'est plus vrai si l'on fait seulement l'hypothèse $v \in \Mt^1_{\alg} $. Par exemple pour $\displaystyle{v=-\frac{z}{2}\d_z}$ et $f=z^2$, on trouve~:
 $$e^v f=z^2-z^2+\frac{z^2}{2!}+\dots+(-1)^n\frac{z^2}{n!}+\dots=e^{-1} z^2$$
alors que $f+v ( f)=0$.  
  
%%%%%%%%%%%%%%%%%%%%
%%%%%%%%%%%%%%%%%%%%
 %%%%%%%%%%%%
 \section{Détermination finie}
 \subsection{Le théorème de Mather-Tougeron}
Considérons  l'anneau local $\Ot_{\CM^n,0}$ des germes de fonctions holomorphes en l'origine dans $\CM^n$.
 Le groupe des automorphismes $\Aut(\Ot_{\CM^n,0})$ de cet anneau agit naturellement. Cet action induit une action infinitésimale
 de l'algèbre de Lie des dérivations $\Der_{\Ot_{\CM^n,0}}(\Ot_{\CM^n,0})$.
 
 Par ailleurs, l'anneau $\Ot_{\CM^n,0}$ est filtré par les puissances de son idéal maximal
 $$\Mt_{\CM^n,0}=\{ f \in \Ot_{\CM^n,0}: f(0)=0 \}. $$
 Le théorème de Mather-Tougeron affirme que pour un germe $f$ dont l'origine est un point critique isolée, pour $k$ suffisamment élevé l'espace affine $f+\Mt^k_{\CM^n,0}$ est contenu dans l'orbite de $f$.
 
 Pour préciser la valeur du nombre $k$, rappelons la définition du nombre de Milnor d'un germe $f \in \Ot_{\CM^n,0}$~\cite{Milnor_isolated}.
 L'{\em idéal jacobien du germe $f$}, noté $Jf$, est l'idéal de $\Ot_{\CM^n,0}$ engendré par les dérivées partielles de $f$. Il revient au même de dire que l'origine est un point critique isolé de $f$ ou bien que l'espace vectoriel
 $\Ot_{\CM^n,0}/Jf$ est de dimension finie. Le nombre 
 $$ \mu(f):=\dim_{\CM} \Ot_{\CM^n,0}/Jf$$
 est appelé le {\em nombre de Milnor} de $f$. Par exemple, pour $n=1$ et 
 $$f:(\CM,0) \to (\CM,0),\ x \mapsto x^2 ,$$
 on $Jf=\Mt_{\CM,0}$, donc  $\mu(f)=1$.
 
 \begin{theoreme}[\cite{Mather_fdet,Tougeron}] Pour tout germe de fonction holomorphe $f:(\CM^n,0) \to (\CM,0)$ dont l'origine est un point critique isolé, l'espace affine $f+\Mt^{\mu(f)+2}_{\CM^n,0}$ est contenu dans l'orbite de $f$ sous l'action du groupe des automorphismes de l'anneau $\Ot_{\CM^n,0}$.
 \end{theoreme}
 \begin{exemple} L'orbite du germe 
 $$f:(\CM,0) \to (\CM,0),\ x \mapsto x^2 ,$$
contient donc tout germe de la forme 
$$(\CM,0) \to (\CM,0),\ x \mapsto x^2+x^3h(x),\ h \in \Mt^3,$$
car $\mu(f)=1$. Plus généralement, tout germe dont la partie quadratique est non-dégénéré peut être ramené à celle-ci par un changement de variable. C'est le lemme de Morse dans un contexte holomorphe.
\end{exemple}
 %%%%%%%%%%%%%%%%%%%%
 \subsection{Le théorème de Poincaré}
 \label{SS::Siegel}
  Le groupe des automorphismes de $\Ot_{\CM^n,0}$ agit sur lui même par conjugaison. Cette action donne lieu
 à la représentation adjointe du groupe des automorphismes dans l'algèbre de Lie des dérivations
 $$\alg:= \Der_{\Ot_{\CM^n,0}}(\Ot_{\CM^n,0}).$$
 Dans les termes classiques~: «tout champ de vecteur peut-être transporté par une application biholomorphe.»
 L'action infinitésimale associée est donné par la représentation adjointe de $\alg$
 $$\alg \to \alg,\ w \mapsto [v,w]. $$
 (On trouve $[v,w]$ et non $[w,v]$ car les automorphismes se composent de gauche à droite contrairement aux applications qui se composent de droite à gauche.)
 
 Notons $z_1,\dots,z_n$ les coordonnées sur $\CM^n$. Le module des dérivations de l'anneau  $\Ot_{\CM^n,0}$
 est librement engendré par les $\d_{z_i}$. La filtration 
 de $\Ot_{\CM^n,0}$ par les puissances de l'idéal maximal donne donc lieu à une filtration de $\alg$~:  
 $$\alg \supset \Mt^1_{\alg} \supset \Mt^2_{\alg} \supset \cdots  .$$
  Le théorème suivant est dû à Poincaré.
\begin{theoreme}[\cite{Poincare_these}]Soit $v \in \Der_{\Ot_{\CM^n,0}}(\Ot_{\CM^n,0})$ une dérivation dont la partie linéaire est diagonalisable. Si  l'enveloppe convexe des valeurs propres ne contient pas l'origine et si celles-ci engendrent un espace vectoriel de dimension $n$ sur $\QM$ alors l'espace affine $v+\Mt_{\alg}^2$ est contenu dans l'orbite de $v$ sous l'action adjointe du groupe $\Aut(\Ot_{\CM^n,0})$.
 \end{theoreme}
%%%%%%%%%%%%%%%%%%%%%%%%%%%%%%%%%
\subsection{Le théorème de Siegel} 
Considérons à présent  le cas où les valeurs propres sont de module égal à un.

 Notons $(\cdot,\cdot)$ la forme bilinéaire dans $\CM^n$ définie par
$$\CM^n \times \CM^n \to \CM,\ (z_1,\dots,z_n,w_1,\dots,w_n) \mapsto \sum_{i=1}^n z_i w_i $$
Pour $i \in \ZM^n$, on pose
$$ \s(i):=| i_1| + |i_2|+ \dots +| i_n |.  $$
Nous dirons qu'un vecteur $\l=(\l_1,\dots,\l_n) \in \CM^n$ est {\em $(C,\tau)$-diophantien},\ $\tau \in \NM$, si
$$\forall i \in \ZM^n \setminus \{ 0 \},\ \left| (\l,i) \right| \geq  \frac{C}{ \s(i)^\tau} $$
et qu'il est {\em diophantien} s'il existe de tels nombres $C,\tau$.

La multiplication par un nombre complexe non-nul envoie l'ensemble des vecteur diophantiens dans lui-même, on peut donc également parler de points diophantiens de l'espace projectif $\PM^{n-1}$. 

Pour $n=2$ et $\a$ algébrique non rationnel, le théorème de Liouville entraîne que le vecteur $\omega=(1,\a)$ est diophantien~\cite{Liouville_approximation1,Liouville_approximation2}. 
Un vecteur $(C,\tau)$-diophantien est un vecteur qui reste «loin» du réseau des entiers.

On vérifie sans difficulté que pour $C,\tau$ fixés, l'ensemble des nombres diophantiens est un ensemble de mesure non-nulle pourvu que 
$\tau>n-1$. De plus, la réunion
de tels ensembles pour les différentes valeurs de $C$, à $\tau$ fixé, forme un ensemble de mesure pleine (voir par exemple \cite{Arnold_edo}).
  
 (Rappelons l'argument~: pour $C,i$ fixés, l'ensemble des nombres qui ne sont pas diophantiens dans $[-N,N]^{2n}$ est un cylindre, la somme sur $i$
 des mesures de ces cylindres converge vers un nombre qui tend vers $0$ avec $C$,  lorsque $\tau>n-1$.)
 
 \begin{theoreme}[\cite{Arnold_SD1,Siegel_linearisation}] Soit $v \in \Der_{\Ot_{\CM^n,0}}(\Ot_{\CM^n,0})$ une dérivation dont la partie linéaire est diagonalisable. Si  les valeurs propres sont diophantiennes alors l'espace affine $v+\Mt_{\alg}^2$ est contenu dans l'orbite de $v$ sous l'action adjointe du groupe $\Aut(\Ot_{\CM^n,0})$.
 \end{theoreme}  
 %%%%%%%%%%%%%
 \subsection{Filtration d'un espace vectoriel échelonné}
 Soit $E$ un espace vectoriel échelonné. Les sous-espaces vectoriels
 $$\Mt^k_E=\{ x \in E: \exists C,\tau,\ | x|_s \leq Cs^k,\ \forall s \leq \tau \}$$
 filtrent l'espace $E$~:
 $$E:=\Mt^0_E \supset \Mt^1_E\supset \Mt^2_E \supset \cdots. $$
 (Dans le cas où l'espace vectoriel $E:=\Ot_{\CM^n,0}$ est échelonné par les espaces de Banach~:
 $$E_s=B(D_s^n),\ | f |_s:=\sup_{z \in D_s^n} | f(z) |,\ D_s^n:=\underbrace{D_s \times \cdots D_s}_{n\ fois}, $$
 on retrouve la filtration de $\Ot_{\CM^n,0}$ par les puissances de son idéal maximal.)
 
 Soit à présent $\alg$ un sous-espace vectoriel de $\Bt^1(E)$ et $G$ un sous groupe de $\Lt(E)$ agissant sur $E$ avec
 $$\exp(\alg) \subset G. $$
 (L'action de $G$ sur $E$ n'est pas nécéssairement l'action naturelle induite par celle de $\Lt(E)$ sur $E$.)
 
 \begin{proposition}
 \label{P::induction}
  Supposons que le groupe $G$ préserve un sous-espace affine $a+M,\ M \subset E$ de telle sorte que
 \begin{enumerate}[{\rm i)}]
 \item $\alg=\Mt^2_{\alg}$ ;
 \item $M$ soit contenu dans la $\alg$-orbite de $a$
 \end{enumerate}
 alors pour tout élément $x \in a+M$ et pour tout $N >0$, il existe $g \in G$ tel que
 $$x=g \cdot a+ z,\ z \in (\Mt_E^N \cap M).$$
 \end{proposition}
 Autrement dit à une correction arbitrairement petite près, tout élément de $a+M$ est contenu dans la 
 $G$-orbite de $a$.
 
 La démonstration est immédiate. \'Ecrivons $x \in E$ sous la forme
 $$x=a+u_0(a). $$
 La première hypothèse entraîne que $u_0$ est exponentiable. En développant l'exponentielle, on trouve
 $$z_1:=x-e^{u_0} a =\sum_{n \geq 0}\frac{u^{n+2}(a)}{(n+2)!}.$$
 Comme $u \in \Mt^2_{\alg}$, on a  
 $$\sum_{n \geq 0}\frac{u^{n+2}}{(n+2)!} \in \alg^4 $$
 et par conséquent $z_1 \in (\Mt_E^3 \cap M)$. En écrivant, $z_1$ sous la forme
 $$z_1=u_1(a) $$
 on a $u_1 \in \alg^8$. On définit  $z_2 \in (\Mt_E^7 \cap M)$ en posant
 $$ z_2:=a-e^{u_1}e^{u_0} a.$$
 On construit ainsi un suite de morphismes bornés $u_0,\dots,u_n$ telle que 
 $$ z_n:=a-e^{u_n} \cdots e^{u_1}e^{u_0} a \in \Mt_E^{2^{n+1}-1}.$$
 Ce qui démontre la proposition.
 
 %%%%%%%%%%%%%%%%%%%
 %%%%%%%%%%%%%%%%%%%%%%%%%%%%%%%%%
\subsection{Comparaison des différents types d'algorithmes itératifs}
Une certaine confusion règne au sujet des différents type d'algorithme itératifs, aussi il n'est peut-être pas tout à fait inutile d'en rappeler les principes.
 
L'algorithme que nous avons utilisé dans la démonstration du n° précédent apparaît dans la démonstration du théorème des tores invariants suggérée par Kolmogrov, nous le nommerons donc {\em algorithme de Kolmogorov}. Il est parfois appelé algorithme de Newton, algorithme de Newton modifié ou encore algorithme KAM.
 
 Pour cela, plaçons nous
dans le cas simple où $G$ est un groupe de Lie, $E$ un espace vectoriel de dimension finie, ici comme par la suite
 $0_E$ désigne l'origine de l'espace vectoriel $E$. 
  
L'action de $G$ sur $E$ définit une action de $T_eG=\alg$ sur $T_{0_E}E \approx E$.  
 On obtient une application 
 $$\rho:\alg \to E,\  \xi \mapsto \xi \cdot 0_E $$
 Supposons que cette application admette un inverse à droite $j$.
Dans l'algorithme de Kolmogorov, on  considère les suites
$$x_0=x,\ u_0=j(x); $$
$$x_n=e^{-u_{n-1}}(x_{n-1}), u_n=j(x_n).$$
Le premier terme de la suite $(x_n)$ vérifie 
$$x_1=e^{-u_0}(x_0)=x_0-u_0 \cdot x_0+({\rm termes\ d'ordre\ sup\acute{e}rieur})  .$$ 
et $x_0-u_0 \cdot x_0=0$ par définition de $u_0$. On recommence cette procédure, en posant $u_1=j(x_1)$, on a alors
$$x_2= e^{-u_1}e^{-u_0}(x_0), u_2=j(x_2)$$
et ainsi de suite...  

Si $(u_n)$ tend vers $0_{\alg}$ et si la suite formée par les produits
$$g_n:= e^{u_n} \dots e^{u_1}  e^{u_0} $$
converge vers une limite $g$ alors la suite $(x_n)=(g_n x)$ converge vers $0_E$. En effet, en passant à la limite
dans l'égalité $u_n=j(x_n) $, on trouve $0_{\alg}=j(gx) $, et par définition de $j$ cela signifie que
$0_E=gx $. Ce qui montre que l'action est transitive.

Remarquons que l'on peut écrire l'itération avec seulement  l'une des deux suites. Pour cela, on considère l'application 
$$\phi :  \alg \rightarrow \alg,
\quad u \mapsto j   \left( e^{-u} u (0_E ) \right).$$
Les itérés  de $u_0=j(x)$ par $\phi$ sont égaux aux termes de la suite $(u_n)$, en effet~:
$$u_{n+1}=j(x_{n+1})=j(e^{-u_n}x_{n})=j(e^{-u_n}u_n ( 0_E)).$$

Dans \cite{groupes}, avec Féjoz, nous avons mis en oeuvre cet algorithme pour les espaces vectoriels échelonnés, mais sous des hypothèses peu commodes et qui ne sont, en règle générale, pas satisfaites pour le paramétrage donné par l'exponentielle.

L'idée de faire intervenir une infinité d'éléments du groupe posait visiblement des problèmes.  Zehnder, tout comme Sergeraert, lui a préféré un algorithme du type Newton. Cet algorithme n'utilise la structure de groupe qu'à posteriori. On a deux espaces vectoriels $E'$ et $E$, on suppose
 qu'une application $f:E' \to E$ admet une différentielle surjective dans un voisinage de l'origine dans $E'$. On fixe $x \in E$ on cherche $u \in E'$ tel que
 $f(u)=x$. On pose 
 $$\p:E' \to E,\ u \mapsto f(u)-x $$
 
 On note $L(u):E \to E'$ un inverse à droite de la différentielle $D\p(u)=Df(u):E \to E'$ de $\p$ au point $u \in E'$. L'algorithme de Newton est donné par les suites
 $$x_0=x,\ u_0=-L(0)x; $$
 $$x_n=\p(u_{n-1}), u_n=u_{n-1}-L(u_{n-1})(x_{n-1})$$
 Si les suites $(u_n)$ et $(x_n)$ convergent respectivement vers $u$ et $x'=\p(u)$ alors $a=a-L(u)(x')$ donc $L(u)x'=0_{E'}$.
 En composant à gauche par $D\p(u)$, les deux membres de cet égalité, on obtient $x'=0$.
 
 Lorsque l'on applique l'algorithme de Newton dans le cas d'une action de groupe
 $$\rho: G  \to E,g \mapsto g \cdot 0_E.$$   la différentielle de $\rho$ au point $g$ est conjuguée à celle au point $e$ où $e$ désigne l'élément unité de $G$. Il suffit donc de connaître les inverses de $\p$ en de tels points.
 
Malheureusement, dans un cadre plus général, cette opération de conjugaison agit sur les estimations nécessaires pour démontrer la convergence de l'itération. Dans certains cas, on peut montrer que cette opération est sans danger \cite[Corollaire 4.2.6]{Sergeraert}. Mais ce sont des cas très exceptionnels. 

Poursuivant une idée de Moser, Zehnder a affaibli la condition d'existence d'un inverse, en la remplaçant par celle d'un inverse approché. Zehnder a tenté de démontrer que pour les action de groupes l'inverse en l'origine était un inverse approché de l'inverse en un point arbitraire, mais il n'a abouti à aucun résultat rigoureux~\cite[Chapter 5]{Zehnder_implicit}.   
 
 On pourrait être tenté de remplacer l'algorithme de Newton par celui de Picard
 $$x_0=x,\ a_0=L(0)x; $$
 $$x_n=\p(a_{n-1}), a_n=a_{n-1}-L(0)(x_{n-1}),$$
 qui ne demande un inverse qu'en l'origine, mais dans ce cas la convergence est plus lente et donc plus hypothétique en dimension infinie.

En résumé, l'algorithme de Kolmogorov est spécifique aux actions de  groupes. Il présente, dans ce cas, les avantages de celui de Newton et de celui de Picard réunis~: c'est un algorithme à convergence rapide qui se construit à l'aide de l'action linéarisée en l'origine.
 
%%%%%%%%%%%%%%%%%%%%%%%%%%%%%%%%%%%%%%%%%%%%%%%%%%%%
\section{Le théorème de $M$-détermination}
 Dans \cite{Zehnder_implicit}, Zehnder écrit~: «(\dots) while all previous proofs of the theorems in question involve infinitely many coordinate changes and consequently complicated convergence arguments, we shall avoid this inconvenience and work in function spaces over a fixed set of variables.» 

 Nous allons donner un critère très simple de convergence qui évite les difficultés évoquées par Zehnder. Celui-ci nous conduira directement à des généralisation des théorèmes de détermination finie et de Poincaré-Siegel.
 %%%%%%%%%%%%%%%%%%%%%%%%%%%%%%%%%%%%%%
\subsection{Produits infinis}
\label{SS::exp}

Soit $E$ un espace vectoriel échelonné. Comme nous allons le voir, le problème de la convergence de l'algorithme de Kolmogorov est résolu par le résultat suivant.
\begin{theoreme} 
\label{T::produits}
Soit $(u_n)$une suite de $\tau$-morphismes $1$-bornés vérifiant la condition (E) du \ref{SS::produits}.
Si la série numérique
$$\sum_{n \geq 0} \frac{N_s^1(u_n)}{s} $$
est convergente pour tout $s \leq \tau$ alors la suite $(g_n)$ définie par
$$g_n:=e^{u_n}e^{u_{n-1}}\cdots e^{u_0}  $$
converge vers un élément inversible de $\Lt(E)$.
 \end{theoreme}
 Pour démontrer ce résultat, commençons par généraliser le corollaire du \ref{SS::produits}.
 \begin{lemme} 
\label{L::produits}
Soit $(u_n)$une suite de $\tau$-morphismes $1$-bornés vérifiant la condition (E).
Pour tout $s \leq \tau$, la norme du morphisme
$$g_n:=e^{u_n}e^{u_{n-1}}\cdots e^{u_0}  $$
vérifie l'inégalité
$$| g_n x |_{\l s} \leq \left( \prod_{i=0}^n \frac{1}{1-\frac{3}{(1-\l)s}N_s^1(u_i)}\right) | x |_s $$
pourvu que $\l$ vérifie
$$\max_{i \leq n}\frac{3}{(1-\l)s}N_s^1(u_i)<1. $$
 \end{lemme}
\begin{proof}
Notons $\D_j \subset \ZM^j$ les suites $ i=(i_1,\dots,i_j)$ dont les éléments sont dans l'ensemble
 $\{0,\dots,n \}$  et telles que $i_p \geq i_{p+1}$.  On a alors la formule
 $$\prod_{i=0}^n \frac{1}{1-z_i}=\sum_{j \geq 0} \sum_{i \in \D_j}z^i,\ z^i:=z_1^{i_1} z_2^{i_2} \dots  z_j^{i_j}$$
 et plus généralement
 $$\prod_{i=0}^n \frac{1}{1-\a z_i}=\sum_{j \geq 0} \sum_{i \in \D_j}\a^ j z^i.$$
 
On pose
$$u[i]:=u_{i_1} u_{i_2} \cdots u_{i_j},\ i \in \D_j. $$
Développons $g_n$ en série puis regroupons les termes suivant l'ordre en $t$, il vient~:
$$g_n=\sum_{j \geq 0} (\sum_{i \in \D_j} u[i]) \frac{t^j}{j!}=1+(\sum_{i=0}^n u_i)t+(\sum_{i=0}^n u_i^2+\sum_{j=0}^n \sum_{i=j+1}^n u_i u_j) \frac{t^2}{2}+\dots .$$
Posons 
$$z_{i,s}:=N_s^1(u_{i})\ {\rm et} \ N_s^j(u[i]):=N_{\tau}^j(u_{i_1} u_{i_2} \cdots u_{i_j}).$$
De la proposition \ref{P::morphismes}, on déduit l'inégalité~:
$$\frac{1}{j!}N_s^j( u[i])  \leq 3^j\prod_{p=0}^jN_s^1(u_{i_p})= 3^j z^i_s,\ i \in \D_j$$
et par suite
$$\left| u[i] x \right|_{\l s} \leq   \left(\frac{ 3}{(1-\l)s} \right)^j z^i_s  | x|_s,\ \forall \l \in ]0,1[. $$
On obtient bien
$$| g_n x |_{\l s} \leq  \left( \sum_{j \geq 0}\sum_{i \in \D_j} \a^j z^i_s \right)| x|_s,\ \a=\frac{3}{(1-\l)s} .  $$
Ce qui démontre le lemme.
 \end{proof}
 
 Achevons la démonstration de la proposition. Pour cela, fixons $s \in ]0,\tau]$ et choisissons $\l \in]0,1[$ tel que
 $$ \sup_{n \geq 0}\frac{3}{(1-\l)s}N_s^1(u_n)) <1.$$
 
 Il est possible de trouver de tels $\l$ car la suite 
 $$\frac{N_s^1(u_n))}{s} $$
 tend vers $0$ et les $u_n$ vérifient la condition (E).

 Montrons tout d'abord que la suite $(g_n)$ définit, par restriction, une suite uniformément bornée d'opérateurs dans $\Lt(\Et_s,\Et_{\l s})$. 
Pour cela, notons $\| \cdot \|_{\l}$ la norme d'opérateur dans $\Lt(\Et_s,\Et_{\l s})$. Le lemme précédent donne l'estimation
$$ \| g_n \|_\l \leq   \prod_{i=0}^n \frac{1}{1-\frac{3}{(1-\l)s}N_s^1(u_i)}.$$

En prenant le logarithme du membre de droite, on voit que le produit converge quand $n$ tend vers l'infini, car la série de terme général
$$\frac{N_s^1(u_n))}{s} $$
est convergente. Comme chacun des facteurs de ce produit est au moins égal à un, on obtient l'inégalité
$$  \| g_n \|_\l \leq  C_\l,\ C_\l:= \prod_{i \geq 0} \frac{1}{1-\frac{3}{(1-\l)s}N_s^1(u_i)}.$$
Ce qui démontre l'assertion.

Soit à présent, $\mu \in ]0,1[$ vérifiant l'inégalité
 $$ \sup_{n \geq 0}\frac{3}{(1-\mu)\l s}N_{\l s}^1(u_n)) <1.$$ Nous allons montrer que la suite $(g_n)$ définit, par restriction, une suite de Cauchy dans $\Lt(\Et_s,\Et_{\mu \l s})$, le théorème en découlera. 

Je dis que la série de terme général $ \| g_n-g_{n-1}\|_{\l \mu }$ est convergente. Pour le voir, écrivons
$$g_n-g_{n-1}=(e^{t u_n}-\Id)g_{n-1} $$
où $\Id \in \Lt(E)$ désigne l'application identité.

En développant l'exponentielle en série, on obtient l'inégalité~:
$$ | (e^{ u_n}-\Id) y |_{\l \mu s} \leq  \left( \sum_{j \geq 0} \frac{ ( 3 N_{\l s}^1(u_n))^{j+1} }{((1-\mu) \l s )^{j+1} } \right) | y |_{\l s}=
\frac{3}{1-\mu-\frac{3N^1_{\l s}(u_N)}{\l s}} \frac{N^1_{\l s}(u_n)}{\l s} | y |_{\l s} ,$$
pour tout $y \in E_{\l s}$. En prenant $y=g_n x$, ceci nous donne l'estimation
$$ \| (e^{ u_n}-\Id) g_{n-1} \|_{\l \mu} \leq   \frac{3 C_\l}{1-\mu-\frac{3N^1_{\l s}(u_n)}{\l s}} \frac{N^1_{\l s}(u_n)}{\l s}.$$
 
La quantité 
$$K_{\l,\mu}:=\sup_{n \geq 0} \frac{3 C_\l}{1-\mu-\frac{3N^1_{\l s}(u_n)}{\l s}} $$
est finie car la suite $\displaystyle{\frac{3N^1_{\l s}(u_n)}{\l s}}$ tend vers $0$ lorsque $n$ tend vers l'infini. 
Nous avons donc montré l'estimation
$$\| g_n-g_{n-1} \|_{\l \mu} \leq K_{\l,\mu} \frac{N^1_{\l s}(u_n)}{\l s}.$$
Il ne nous reste plus qu'à utiliser l'inégalité triangulaire pour voir que $(g_n)$ définit une suite de Cauchy de l'espace de Banach $\Lt(\Et_s,\Et_{\mu \l s})$~:
$$\| g_{n+p}-g_{n}  \|_{\l \mu } \leq \sum_{i = 1}^{p}  \| g_{n+i}-g_{n+i-1}\|_{\l \mu  } \leq K_{\l,\mu}\left(\sum_{i = 1}^{p} \frac{3N^1_{\l s}(u_{n+i})}{\l s}\right).$$
Nous avons donc montré que la suite $(g_n)$ converge vers un élément $g \in \Lt(E)$. On démontre de même que la suite $(h_n)$  définie par
$$h_n=e^{- u_0}e^{- u_1} \cdots e^{-u_n}$$ 
converge vers un élément $h \in \Lt(E)$. Pour tout $n \in \NM$, on a
$$g_nh_n=h_ng_n=\Id$$
 donc $gh=hg=\Id$. Ce qui montre que $h$ est l'inverse de $g$.
 Le théorème est démontré.
 %%%%%%%%%%%%%%%%%%%
 \subsection{\'Enoncé du théorème de $M$-détermination}
 Soit $E$ un espace vectoriel $S$-échelonné et $M$ un sous-espace vectoriel de $E$.
 \begin{definition}L'action d'un sous-espace vectoriel $\alg \subset \Bt^1(E)$ est dite échelonnée en $a \in E$ s'il existe un constante
 $C >0$ telle que
 $$| u^n \cdot a|_s\leq \frac{C}{\s^n}N_{s+\s}^n(u^n) $$
 pour tous $u \in \alg$, $n \geq 1$, $\s \in ]0,S-s[$.
 \end{definition}
 \begin{exemple} Soit $E=(E_s)$ un espace vectoriel $S$-échelonné et $a \in E_{S'}$. Considérons le nouvel échelonnement $(E'_s)$ de $E$ obtenu
 en remplaçant l'intervalle d'échelonnement $[0,S]$ par $[0,S']$~: 
 $$E_s':=E_s,\ s \leq S'.$$
 L'action naturelle de $\Bt^1(E)$ est échelonnée, il suffit de prendre 
 $$C:=| a |_{S'}. $$
 Par conséquent, quitte à modifier l'intervalle d'échelonnement, on peut toujours supposer que l'action naturelle des morphismes bornés est échelonnée en un point donné.
 \end{exemple}
 \begin{theoreme}
 \label{T::groupes} Soit $E$ un espace vectoriel échelonné, $a \in E$, $M$ un sous-espace vectoriel fermé de $E$, $\alg$ un sous-espace vectoriel de $\Bt^1(E)$, $G$ un sous-groupe fermé de $\Lt(E)$ contenant $\exp(\alg)$ et agissant sur $a+M$. Si
 \begin{enumerate}[{ \rm i)}]
 \item l'action de $\alg$ est échelonnée en $a$ ;
 \item l'application  $\rho:\alg \to M,\ u \mapsto u \cdot a$ possède un inverse borné ;
 \item $\alg=\Mt^2_{\alg}$
 \end{enumerate}
 alors l'orbite de $a$ sous l'action de $G$ est égale à $a+M$.
 \end{theoreme}
%%%%%%%%%%%
 \subsection{Principe de la démonstration du théorème \ref{T::groupes}}
 Notons 
 $$j:E \mapsto \alg$$
   l'inverse de l'application $\rho$. Soit $b \in M$, on cherche $g \in G$ tel que  $g \cdot a=a+b$.
Pour cela, on  considère les suites $(b_n)$ et $(u_n)$ définies par~:
\begin{enumerate}[1)]
\item $b_{n+1}:=e^{-u_n}(a+b_n)-a$ ;
\item  $u_{n+1}:=j(b_{n+1}).$
\end{enumerate}
avec $ b_0=b,\ u_0=j(b) $.

Dans cette itération, la suite $(u_n)$ peut également être définie par la formule
$$u_{n+1}=j((e^{-u_n}(a+u_n\, a )-a)) $$
 
 Supposons que la suite formée par les produits
$$g_n:= e^{u_n} \dots e^{u_1}  e^{u_0} $$
converge vers une limite $g'$ et que $(u_n)$ tende vers $0_{\alg}$. Dans ce cas, la suite $(x_n)=(g_n x)$ converge vers $x'=a+b'$. Par définition de $j$, on a
$$ u_n(a)=b_n$$
et en passant à la limite sur $n$, dans les deux membres de l'égalité, on trouve 
$$0_E=b'.$$ 
Ce qui montre que $g' \cdot (a+b)=a$ donc l'élément $g$ que nous cherchions est l'inverse de $g$. CQFD.
 
Le théorème sera donc démontré pourvu que $(g_n)$ soit convergente et que $(u_n)$ tende vers $0_{\alg}$.  D'après le théorème \ref{T::produits},  il suffit pour cela de  montrer que la série
$$ \sum_{n \geq 0}\frac{N^1_s(u_n)}{s}$$ est  convergente pour tout $s$ suffisamment petit, et comme nous allons le voir c'est un fait presque immédiat. 

%%%%%%%%%%%%%%%%%%%%%%%%%%%%%%%%%
\subsection{Démonstration du théorème \ref{T::groupes}}
\label{SS::homogene_dem}

Comme l'action de $\alg$ est échelonnée en $a \in E$, il existe une constante $C>0$ telle que
$$| u^n \cdot a |_s \leq \frac{C}{\s^n}N_{s+\s}^n(u^n) $$
pour tout $u \in \alg_\tau$ et tout $\s \in ]0,\tau-s[$. 
 
 \begin{lemme}
 \label{L::reste}
 Pour tout $\tau$-morphisme $1$-borné $u$ vérifiant la condition 
 $$\frac{3N^1_\tau(u)}{\tau-s} \leq \frac{1}{2}, $$
  on a l'inégalité~:
 $$| (e^{-u}(\Id +u)-\Id) \cdot a|_s \leq  \frac{36C }{(\tau-s)^2} N^1_\tau(u)^2,$$
 pour tout $s \in ]0,\tau[$.
\end{lemme}  
\begin{proof}
On a l'égalité~:
$$ e^{-u}(\Id +u)-\Id=\sum_{n \geq 0} \frac{(n+1)}{(n+2)!}(-1)^{n+1}u^{n+2}.$$
Comme l'action de $\alg$ est échelonnée, on en déduit l'estimation 
$$ | \sum_{n \geq 0} (-1)^{n+1} \frac{(n+1)}{(n+2)!}u^{n+2} \cdot a |_s \leq C \sum_{n \geq 0}  
\frac{(n+1)3^{n+2}}{(\tau-s)^{n+2}}N^1_\tau(u)^{n+2} .$$
Comme
$$\frac{3N^1_\tau(u)}{\tau-s} \leq 1 $$
le membre de droite est égal à
$$Cx^2\sum_{n \geq 0}(n+1)x^n=\frac{Cx^2}{(1-x)^2},\ {\rm \ avec\ }x= \frac{3N^1_\tau(u)}{\tau-s}.$$
En utilisant l'inégalité
$$\frac{1}{(1-x)^2} \leq 4,\ \forall x \in [0,\frac{1}{2}], $$
on trouve bien la majoration du lemme.
 \end{proof}

 Soit $k$ tel que $j$ soit $k$-borné. Quitte à réduire l'intervalle d'échelonnement, on peut supposer que $j$ est un $S$-morphisme $k$-borné avec
 $S<1$.
  Comme $u_0$ est d'ordre $2$, le lemme précédent montre, en particulier, que $u_k$ est d'ordre $2^{k+1}$. Par conséquent,  tout $s$ suffisament petit, vérifie l'inégalité~:
$$(*)\ N^1_{2s}(u_{k+1}) \leq  2^{-3(k+2)}m \left(\frac{s}{8} \right)^{k+2},\ m:=\min(1,\frac{1}{36CN^k_S(j)}).  $$
Nous allons montrer que si $s$ est de la sorte et si $a \in E_{2s}$ alors la série 
$$ \sum_{n \geq 0}\frac{N^1_s(u_n)}{s}$$ est  convergente.

Pour cela, considérons la suite $(\s_n)$ définie par $\s_n=0$ pour $n \leq k$ et
$\displaystyle{\s_n=\frac{s}{2^{n-k+1}},}$ pour $n \geq k+1$~:
$$ \s_0=\dots=\s_k=0,\ \s_{k+1}=\frac{s}{4},\  \s_{k+2}=\frac{s}{8},\dots $$
Définissons la suite $(s_n)$ par 
$$ s_{n+1}=s_n-2\s_n,\ s_0=2s .$$
On a
$$s_0=s_1\cdots=s_{k+1}=2s,\ s_{k+2}=\frac{3}{2}s,\ s_{k+3}=\frac{5}{4}s,\cdots $$
 
 Montrons par récurrence sur $n$ que, pour $n \geq k+1$, on a~:
$$(**)\ N^1_{s_n}(u_n) \leq m\s_{n+1}^{k+2}=2^{-(n-k+2)(k+2)}m s^{k+2} . $$
Cette inégalité est vérifié pour $n=k+1$ d'après (*), il nous reste donc à montrer que l'inégalité au rang $n \geq k+1$ implique l'inégalité au rang $(n+1)$.

Appliquons le lemme avec
$$b_{n+1}:=(e^{-u_n}(\Id +u_n)-\Id)\cdot a $$
et $\tau-s=\s_n$. On obtient l'inégalité~:
$$| b_{n+1}|_{s_n-\s_n} \leq  \frac{36C }{\s_n^2} N^1_{s_n}(u_n)^2.$$
En utilisant l'hypothèse de récurrence et la définition de $m$, on obtient l'estimation
$$| b_{n+1}|_{s_n-\s_n} \leq \frac{m \s_{n+1}^{2k+4}}{N^k(j)\s_n^2}=  \frac{m \s_{n+1}^{2k+2}}{4N^k(j)}$$
Comme $j$ est $k$-borné, et $u_{n+1}=j(b_{n+1})$, on en déduit l'inégalité
$$N^1_{s_{n+1}}(u_{n+1})\leq \frac{m \s_{n+1}^{2k+2}}{4\s_n^k}= m \s_{n+2}^{k+2}. $$ 
Le théorème est démontré.
%%%%%%%%%%%%%%%%%%%%%%%%%
 \section{Théorème de transversalité}
 \subsection{Transversales}
 Soit $E$ un espace vectoriel, $M$ un sous-espace vectoriel de $E$ et $a \in E$. Soit $G$ un groupe agissant sur $a+M$.
 Nous dirons qu'un sous espace affine $a+F \subset a+M$ est une  {\em transversale} à l'orbite de $a$ sous l'action de $G$ dans $a+M$ si  l'application
  $$G \times F \to a+M,\ (g,\a) \mapsto g(a+\a) $$
   est surjective.   Dans le cas particulier où $F$ est invariant par l'action de $G$, il revient au même de dire que $a+M$ est un espace $(G \times F)$-homogène.

L'espace $F$ est une  {\em transversale} à l'orbite de $a$ sous l'action de $\alg$ si l'application
  $$\alg \times F \to F \times E/F,\ (\a,u) \mapsto (\a,\overline{u(\a)}) $$
  est surjective. Dans ce cas, pour tout $\a \in F$, tout élément $x \in E$ s'écrit sous la forme
$$x=u(\a)+\b,\ \a,\b \in F .$$
 Ici et par la suite, l'application
 $$E \to E/F,\ x \mapsto \overline{x} $$
 désigne la projection canonique.
 %%%%%%%%%%%%%%%%%%%%%
%%%%%%%%%
\subsection{\'Enonc\'e du théorème de transversalité}
Nous désignons par $E$ un espace vectoriel $S$-échelonné et  $F,M$ des sous-espaces vectoriels fermés de $E$ avec $F \subset M$.
 \begin{definition}L'action d'un sous-espace vectoriel $\alg \subset \Bt^1(E)$ est dite $F$-échelonnée en $a \in E$ s'il existe un constante
 $C >0$ telle que~:
 \begin{enumerate}[{\rm A)}]
 \item $\displaystyle{| u^n \cdot a|_s\leq \frac{C}{\s^n}N_{s+\s}^n(u^n) }$ ;
 \item  $\displaystyle{| u^n \cdot \a|_s\leq \frac{C| \a|_{s+\s}}{\s^n}N_{s+\s}^n(u^n) }$
 \end{enumerate}
 pour tous $u \in \alg$, $n \geq 1$, $s \in ]0,S[$, $\s \in ]0,S-s[$, $ \a \in F_s$.
 \end{definition}
 Les actions échelonnés du § précédent correspondent au cas $F=\{ 0 \}$.
 
Nous aurons également besoin de la notion d'application bornée pour les applications qui ne sont pas nécessairement linéaires.
\begin{definition} Une application $f:E \to F$ entre deux espaces vectoriels $S$-échelonnés est dite $k$-bornée si
\begin{enumerate}[{\rm a)}]
\item l'image de $E_s$ par $f$ est contenue dans $F_{s'}$ pour tout $s'<s$ ;
\item il existe une constante $N$ telle que
$$ | f(x) |_s \leq  \frac{N}{\s^k}(1+| x |_{s+\s})$$
pour tous $s \in ]0,S[$, $\s \in ]0,S-s[$.
\end{enumerate}
\end{definition}
Comme pour les applications linéaires nous dirons d'une application qu'elle est bornée s'il existe $k$ pour lequel elle est $k$-borné.
\begin{exemple} Si $P \in \RM[X]$ est un polynôme et $u:E \to F$ un morphisme borné entre espace vectoriels échelonnés alors $P(u)$ est une application affine bornée.
\end{exemple}

\begin{theoreme}
 \label{T::transversale} Soit $M$ un sous-espace vectoriel fermé de $E$, $\alg$ un sous-espace vectoriel de $\Bt^1(E)$, $G$ un sous-groupe fermé de $\Lt(E)$ contenant $\exp(\alg)$ et agissant sur $a+M$ avec $\alg=\Mt^2_{\alg}$ et $F=\Mt^2_F$. Supposons que  l'action de $\alg$ soit $F$-échelonnée en $a$ et qu'il existe pour un certain $k \geq 0$ une application  bornée 
 $$j:F \mapsto \Bt^k(M/F,\alg)$$
  telle que   $j(\a)$ soit un inverse de
$$ \alg \to M/F,  u \mapsto u \cdot (a+\a) $$ 
 alors $F$ est une transversale à l'orbite de $a$ sous l'action de $G$ dans $a+M$.
 \end{theoreme}
 
 %%%%%%%%%%%%%
 \subsection{Principe de la démonstration du théorème \ref{T::transversale}}
 Pour tout $b \in E$, il s'agit de trouver $g \in G$ et $\a \in F$ tels que
 $$g \cdot (a+b)=a+\a. $$
 Pour cela considérons les suites $(a_n)$, $(b_n)$, $(\a_n)$ et $(u_n)$ définies par
 \begin{enumerate}[{\rm i)}]
 \item $a_{n+1}=a_n+\a_n$ ;
 \item $b_{n+1}=e^{-u_n}\cdot (a_n+b_n)-a_{n+1}$ ;
 \item $u_{n+1}=j(\a_n)\cdot b_{n+1}$ ;
  \item $\a_n=b_n-u_n \cdot a_n$ ;
 \end{enumerate}
 avec  $a_0=a$, $b_0=b$, $u_0=j(a)b$. Pour simplifier, les notations nous poserons $\a_{-1}=0$.
 
 Si les $u_i$ sont exponentiables, si la suite  formée par les produits
 $$g_n:=e^{-u_n} e^{-u_{n-1}} \cdots e^{-u_0} $$
 converge vers une limite $g'$, si les suites $(u_n)$ $(b_n)$  tendent vers $0$ et si $(\a_n)$ est sommable alors en passant à la limite dans l'égalité
 $$a_{n+1}+b_{n+1}=g_n \cdot (a+b) $$
 on trouve
 $$g' \cdot (a+b)=a+\sum_{n \geq 0} \a_n \in a+F $$
  
%%%%%%%%%%%%%%%%%%%
\subsection{Démonstration du Théorème \ref{T::transversale}}
Quitte à remplacer l'intervalle d'échelonnement $]0,S[$ par un intervalle plus petit, on peut supposer que $S<1$.

Comme l'action de $\alg$ est $F$-échelonnée en $a \in E$, il existe une constante $C>0$ telle que~:
\begin{enumerate}[{\rm 1)}]
 \item $\displaystyle{| u^n \cdot a|_s\leq \frac{C}{2\s^n}N_{s+\s}^n(u^n) }$ ;
 \item  $\displaystyle{| u^n \cdot \a|_s\leq \frac{C| \a|_{s+\s}}{2\s^n}N_{s+\s}^n(u^n) }$, 
 \end{enumerate}
 pour tous $u \in \alg$, $n \geq 1$, $s \in ]0,S[$, $\s \in ]0,S-s[$ et $\a \in F_{s+\s}$.

Si $|\a|_{s+\s} \leq 1$ alors ces deux conditions donnent lieu à l'inégalité~:
$$| u^n \cdot (a+\a)|_s\leq \frac{C}{\s^n}N_{s+\s}^n(u^n). $$

Pour $n>0$, le vecteur $b_n \in M$ s'écrit sous la forme
$$b_n=A_n+B_n $$
avec
$$A_n:=(e^{-u_{n-1}}(\Id+u_{n-1})-\Id) a_{n-1},\ B_n:=(e^{-u_{n-1}} -\Id) \a_{n-1}. $$

Comme $F=\Mt^2_{F}$ et $\alg=\Mt^2_{\alg}$, pour $i>0$, on a~:
$$u_i \in \Mt^{2^{i+1}}_{\alg} \ {\rm et\ } \a_i,A_i,B_i \in \Mt_E^{2^{i+1}}. $$

Fixons $l$ tel que $j$ soit $l$-borné et  soit $N(j)$  tel que
$$2\s^ lN^k_s(j(\a)) \leq N(j)(1+ | \a |_{s+\s} )$$
pour tous $ s \in ]0,S[$, $\s \in ]0,S-s[$ et $  \a  \in F_{s+\s}$.

Je dis que si $s$ est suffisamment petit pour que les inégalités suivantes aient lieu~:
  \begin{enumerate}[{\rm a)}]
  \item $\displaystyle{| A_{k+l+2}|_{3s} \leq \frac{ m^2 }{N(j)}} \left(\frac{s^{2(k+l)+3}}{9} \right) $ ;
   \item $\displaystyle{| B_{k+l+2}|_{3s} \leq \frac{ m^2}{N(j)}} \left(\frac{s^{2(k+l)+3}}{9} \right) $ ;
 \item  $\displaystyle{| \a_{k+l} |_{3s} \leq \frac{ms^{k+l+1}}{27}}$ ;
 \item  $\displaystyle{| \a_i |_{3s} \leq \frac{1}{3(k+l+1)}}$ pour $i \leq k+l$, 
 \end{enumerate}
 avec
 $$m:=\min(1,\frac{1}{36CN(j)},\frac{1}{2^3.3^5C},\frac{1}{2N(j)}),$$
alors l'algorithme décrit au n° précédent est convergent.

Considérons la suite $(\s_n)$ définie par $\s_n=0$ pour $n \leq k+l,$
 et $\displaystyle{\s_n=\frac{s}{3^{n-k-l}}}$
pour $n \geq k+l+1$~:
$$\s_0=\s_1=\cdots=\s_{k+l}=0,\ \s_{k+l+1}=\frac{s}{3},\ \s_{k+l+1}=\frac{s}{9},\dots $$

Définissons la suite $(s_n)$ par
$$ s_{n+1}=s_n-3\s_n,  s_0=\frac{5s}{2}$$
de telle sorte que $s=(s_n)$ décroît de $\displaystyle{s_n=\frac{5s}{2}}$ vers $s$:
$$s_0=s_1= \dots=s_{k+l+1}=\frac{5s}{2},\ s_{k+l+2}=\frac{3s}{2},\ s_{k+l+3}=\frac{7s}{6},\dots $$

Montrons par récurrence sur $n$ que pour tout $n \geq k+l+2$, on a
 \begin{enumerate}[{\rm i)}]
 \item $\displaystyle{| A_n|_{s_{n-1}-\s_{n-1}} \leq m^{n-k-l}\frac{ \s_n^{2(k+l)+3} }{2N(j)}} $ ;
  \item $\displaystyle{| B_n|_{s_{n-1}-2\s_{n-1}} \leq m^{n-k-l}\frac{ \s_n^{2(k+l)+3} }{2N(j)}} $ ;
 \item $\displaystyle{| \a_n |_{s_n-\s_n} \leq m\s_{n+1}^{k+l+1} }$.
 \end{enumerate}
Supposons i), ii) et iii) vérifiés jusqu'au rang $n$.
On a $\a_{n-1} \in E_{s_n+2\s_{n-1}}$ car $s_{n-1}-\s_{n-1}=s_n+2\s_{n-1}$. Comme $j$ est $l$-borné, on en déduit que~:
$$N_{s_n+2\s_{n-1}}^k (j(\a_{n-1})) \leq \frac{N(j)}{2\s_{n-1}^l}(1+| \a_{n-1} |_{s_n+\s_{n-1}}) \leq \frac{N(j)}{\s_{n-1}^l} $$
 
 Comme $u_n=j(\a_{n-1})(A_n+B_n), $ on a donc~:
 $$N^1_{s_n}(u_n) \leq \frac{N(j)}{\s_{n-1}^{k+l}} | A_n+B_n|_{s_n+\s_{n-1}}$$
  pour tout $\displaystyle{s <\frac{S}{3}}.$ 
 Comme $s_{n-1}-2\s_{n-1}=s_n+\s_{n-1}$, les deux premières inégalités donnent lieu à l'estimation
 $$(*)\ N^1_{s_n}(u_n) \leq  \frac{m^{n-k}\s_n^{2(k+l)+3}}{\s_{n-1}^{k+l}}  = 27m^{n-k}\s_{n+1}^{k+l+3} $$
Par conséquent, le théorème sera démontré pourvu que i), ii) et iii) le soit.

Pour $n=k+l+2$, les inégalités i), ii) et iii) résultent de a), b) et c). Montrons que la validité ces inégalités au rang $n$ entraîne celle au rang $(n+1)$.
En appliquant le lemme du \ref{SS::homogene_dem} avec
$$A_{n+1}:=(e^{-u_n}(\Id +u_n)-\Id) a_n $$
et $\tau-s=\s_n$, on obtient une estimation plus précise que i). En effet, en utilisant b) et ii) jusqu'au rang $n-1$, on a
 $$\sum_{i=0}^{n-1}| \a_i |_{s_n} \leq 1 .$$
 Comme l'action de $\alg$ est $F$-échelonnée,  on obtient ainsi l'inégalité~:
$$| A_{n+1}|_{s_n-\s_n} \leq  \frac{36C}{\s_n^2} N^1_{s_n}(u_n)^2.$$ 
En utilisant l'inégalité (*) au rang $n$, on trouve
$$| A_{n+1}|_{s_n-\s_n} \leq   \frac{2^2.3^8Cm^{2(n-k-l)-1} \s_{n+1}^{2(k+l+3)}}{\s_n^2} =2^2.3^5C m^{2(n-k-l)-1} \s_{n+1}^{2(k+l)+4},   $$
qui est bien, par définition de $m$, une estimation plus précise que i).
  
 Pour montrer la deuxième inégalité, nous aurons besoin du
  \begin{lemme} Pour tout $\a \in F$ et pour tout $\tau$-morphisme $1$-borné $u \in \alg$ vérifiant la condition (E), on a l'inégalité~:
 $$| (e^{-u}-\Id) \cdot \a|_s \leq  \frac{6C |\a|_\tau}{(\tau-s)} N^1_\tau(u),$$
 pour tout $s \in ]0,\tau[$ tel que 
 $$\frac{3N_\tau^1(u)}{\tau-s} \leq \frac{1}{2}. $$
\end{lemme}  
\begin{proof}
L''égalité~:
$$ e^{-u}-\Id=\sum_{n \geq 0} \frac{1}{(n+1)!}(-1)^{n+1}u^{n+1}.$$
Comme l'action de $\alg$ est $F$-échelonnée, on en déduit l'estimation~: 
$$ | (\sum_{n \geq 0} \frac{1}{(n+1)!}(-1)^{n+1}u^{n+1}) \a|_s \leq C |\a|_{\tau} \sum_{n \geq 0}  
\frac{3^{n+1}}{(\tau-s)^{n+1}}N_\tau(u)^{n+1} .$$
Comme
$$\frac{3N_\tau^1(u)}{\tau-s} \leq 1 $$
le membre de droite est fini et égal à
$$\frac{C| \a |_s x}{1-x},\  x= \frac{3N_\tau^1(u)}{\tau-s} $$
En utilisant l'inégalité
$$\frac{1}{(1-x)} \leq 2,\ \forall x \in [0,\frac{1}{2}], $$
on trouve bien la majoration du lemme.
 \end{proof}
 Appliquons ce lemme à $u_n,\a_n$ avec $\tau-s=\s_n$. On obtient l'inégalité~:
$$| B_{n+1}|_{s_n-2\s_n} \leq  \frac{6C}{\s_n} N_{s_n}(u_n) |\a_n|_{s_n-\s_n}.$$
L'hypothèse de récurrence iii) et l'estimation (*) au rang $n$ entraînent l'inégalité~:
$$| B_{n+1}|_{s_n-2\s_n} \leq  2.3^4 Cm^{n-k-l+1}\s_{n+1}^{2(k+l)+3}   . $$
Par définition de $m$, c'est un inégalité plus précise que ii) au rang $(n+1)$. Ces estimations des normes de $A_{n+1}$
et de $B_{n+1}$ entraînent en particulier, l'inégalité (*) au rang $n+1$.
 
Enfin, remarquons que $\a_{n+1}$ est donné par la formule
$$\a_{n+1}=b_{n+1}-u_{n+1}\cdot a_{n+1} $$
et par suite
$$|\a_{n+1}|_{s_{n+1}-\s_{n+1}} \leq   m^{n-k-l+1}\frac{\s_{n+1}^{2(k+l)+3}}{N(j)}+ |u_{n+1} \cdot a_{n+1} |_{s_{n+1}-\s_{n+1}}. $$
Comme l'action de $\alg$ est échelonnée, d'après c) et iii) jusqu'au rang $n$, on a~:
 $$\sum_{i=0}^n| \a_i |_{s_{n+1}} \leq 1 ,$$
 on en déduit
$$ |u_{n+1} \cdot a_{n+1} |_{s_{n+1}-\s_{n+1}} \leq \frac{C}{\s_{n+1}}N^1_{s_{n+1}}(u_{n+1})\leq \frac{3^3Cm^{n-k-l+1} \s_{n+2}^{k+l+3}}{\s_{n+1}}=3^2C   m^{n-k-l+1}\s_{n+2}^{k+l+2} . $$
Ce qui nous conduit également, par définition de $m$, à une estimation plus précise que iii).
Le théorème est démontré.

%%%%%%%%%%%%
\section{Déformations d'espaces vectoriels échelonnés}
 %%%%%%%%%%%%%%%%%
  %%%%%%%%%%%%%%%%%%%%%%%%%ANCIEN
%%%%%%%%%%%%%%%%%%%%%%%%%
 
\subsection{Stabilité des fonctions de Morse}
En théorie des singularités, on distingue les théorèmes de formes normales et ceux de stabilité par déformations. Prenons le cas simple de l'anneau $\Ot_{\CM^n,0}$. Rappelons qu'un germe de fonction
$$f:(\CM^n,0) \to (\CM,0) $$ 
est dit de Morse si l'origine est un point critique dont la partie quadratique est non dégénérée. Pour $n=1$, cela signifie simplement que la série de Taylor de $f$ est de la forme
$$f(x)=ax^2+\dots,\ a \neq 0.$$  
On peut énoncer le lemme de Morse sous trois formes~:
\begin{enumerate}[{\rm i)}]
\item toute germe de fonction  proche d'une fonction de Morse $f$ se ramène à $f$ par un changement de variables ;
\item l'espace affine $f+\Mt^2_{\CM,0}$ est contenu dans l'orbite de $f$ pour toute fonction de Morse $f$ ;
\item toute déformation d'une fonction de Morse est triviale.
\end{enumerate}
(Dans ce contexte, une déformation d'un germe $f$ est un germe d'application de la forme
$$F:(\CM^k\times \CM^n,0) \to (\CM,0),\ F(\l=0,\cdot)=f.  {\rm )}$$

Montrons que les deux premières formes sont équivalentes et que la troisième est plus forte.

\noindent i) $\implies$ ii) est évident et l'équivalence provient du fait que tout germe de fonction proche d'une fonction de Morse possède un point
critique de Morse dans un voisinage de l'origine.

\noindent iii) $\implies$ ii) s'obtient en interpolant linéairement $f+g$ et $f$:
$$f_\e=(1-\e)f+\e g. $$ 
Pour chaque $\e \in [0,1]$, le germe d'application
$$(\CM^2,0) \to (\CM,0),\ (t,x) \mapsto (1-\e-t)f+(t+\e)g $$
est un déformation triviale de $f_\e$. Par conséquent, tout $\e \in [0,1]$ admet un voisinage $V_\e$ dans lequel $f_\e$ et $f_{\e'}$ se ramène à $f_\e$ par un changement de variables, pour $\e' \in V_\e$.  La propriété de Borel-Lebesgue pour l'intervalle $[0,1]$ permet de déduire $f_1$ de $f_0$ par un changement de variable. (Par cet astuce, on évite de montrer que les solutions des équations cohomologiques sont valables pour tout $\e \in[0,1]$.).  

%%%%%%%%%%%%%%%%%
\subsection{Le théorème de stabilité de Mather}
Les germes de transformations biholomorphe de $\CM^n$ et $\CM^m$ agissent par composition sur l'espace des germes
d'applications de $\CM^n$ dans $\CM^m$~:
$$\xymatrix{ (\CM^n,0) \ar[d]^\p \ar[r]^f & (\CM^m,0) \ar[d]^\psi \\
                       (\CM^n,0) \ar[r]^{(\p,\psi)\cdot f} & (\CM^m,0)}$$
On en déduit une action du groupe $\Aut(\Ot_{\CM^n,0}) \times \Aut(\Ot_{\CM^m,0})$ sur l'espace $\Mt_{\CM^n,0}^m$ des germes d'applications holomorphes
$$f:(\CM^n,0) \to (\CM^m,0) $$
qui s'annulent en l'origine.
 
 On peut étendre ces considérations aux déformations des applications. Nous notons $\Aut(\Ot_{  \CM^{n+k}/\CM^k,0})$, $k \leq n$, le groupe des automorphisme de l'anneau $\Ot_{  \CM^{n+k},0}$ qui commutent à la projection sur les $k$premières variables~:
$$\xymatrix{ (\CM^k \times \CM^n,0) \ar^{\p}[r] \ar^{\pi}[d] & (\CM^k \times \CM^n,0) \ar[d]^\pi \\
                     (\CM^k,0) \ar[r]^{\widetilde \p}  & (\CM^k,0) }. $$
Les dérivations correspondantes seront notés $\Der(\Ot_{\CM^{n+k}/\CM^k,0})$.
L'action du groupe $G$ entraîne une {\em action infinitésimale} de 
$$\alg:=\Der(\Ot_{\CM^{n+k}/\CM^k,0}) \times \Der(\Ot_{  \CM^{m+k}/\CM^k,0}).$$
Notons respectivement $z_1,\dots,z_n$ et $w_1,\dots,w_m$ les coordonnées dans $\CM^n$ et dans $\CM^m$. L'action infinitésimale
donnée par la formule
$$ (\sum_{i=1}^n a_i \d_{z_i},\sum_{i=1}^m b_i \d_{w_i}) \mapsto \sum_{i=1}^n a_i \d_{z_i}  f+\sum_{i=1}^m b_i(f_i)  $$

Une déformation d'un germe $f$ est dite {\em triviale} si elle est dans l'orbite de $f$ sous l'action  du groupe $\Aut(\Ot_{  \CM^{n+k}/\CM^k,0})$. Un germe est dit {\em stable} si toutes ses déformations sont triviale. Par exemple, les germes de fonctions de Morse sont stables et pour $m=1$ ce sont les seules fonctions stables ayant un point critique en l'origine. (Bien entendu, le théorème des fonctions implicites montre que les germes dont la dérivée est de rang maximal sont stables.)

\begin{theoreme}[\cite{Mather_stability}] Pour qu'un germe d'application holomorphe $f:(\CM^n,0) \to (\CM^m,0)$ soit stable, il faut et il suffit que
l'orbite de $f$ sous l'action de $\alg$ soit égale à  $\Mt_{\CM^n,0}^m$.
\end{theoreme}
\begin{exemple}
Dans \cite{Whitney_singularities}, Whitney a démontré la stabilité la fronce~:
$$(\CM^2,0) \to (\CM^2,0),\ (z_1,z_2) \mapsto (\frac{z_1^3}{3}+z_2 z_1,z_2). $$
Montrons que cela résulte du théorème de Mather (ici $m=n=2$).
Tout élément $(a,b) \in \Mt_{\CM^n,0}^2 $ possède un élément de la forme $(a',0)$ dans sa $\alg$ orbite. En effet
$$(a,b)-b\d_{z_2}f=(a-b\d_{z_2}f,0). $$
Tout germe de la forme $(a,0)$  possède un élément de la forme $(b(z_2)+c(z_2)z_1,0)$ dans son orbite. Pour le voir écrivons $a$ sous la forme
$$a=a_0(z_2)+a_1(z_2)z_1+a_2(z_1,z_2)z_1^2. $$
L'élément $(a,0)-a_2\d_{z_1}f$ est de la forme anoncée. Finalement
$$(b(z_2)+c(z_2)z_1,0)-c(z_2)\d_{z_2}f=(b(z_2),c(z_2))=b\d_{w_1}f+c\d_{w_2} f .$$
Donc tout élément de $\Mt_{\CM^2,0}^2$ est contenu dans l'orbite de $f$, ce qui démontre l'affirmation.
\end{exemple}

%%%%%%%%%%
 \subsection{Définition, exemple}
 Soit $E$ un espace vectoriel $S$-échelonné.
 Nous dirons d'une paire $(\Et,t)$ formé d'un espace vectoriel $S$-échelonné $\Et$ et d'une application linéaire $t \in \Lt(\Et)$ que c'est une {\em déformation de $E$}   si elle satisfait aux conditions suivantes~:
 
\noindent ($E_1$) L'application linéaire $t \in \Lt(\Et)$ induit une suite exacte
$$0 \to \Et \stackrel{t}{\to} \Et \to E \to 0. $$
($E_2$) L'application linéaire $t$ est un morphisme $0$-borné qui vérifie l'inégalité $N^0_s( t)  \leq s^2$ pour tout $s \in ]0,S[$.

 Dans les cas les plus simples, $\Et$ est une algèbre commutative et l'application linéaire $t$ est la multiplication par un des éléments du centre de l'algèbre. Nous noterons simplement $tx$ l'image de $x \in \Et$ par le morphisme $t \in \Lt(\Et) $

Voici un tel exemple. On note $\Et_s$ l'espace des fonctions continues sur le polycylindre $P_s:=D_s \times D_{s^2}$ et holomorphes dans  son intérieur. On munit $\Et_s$ d'une structure d'espace de Banach en posant
 $$| f |_s:=\sup_{z \in P_s} | f(z_1,z_2) | .$$
 En prenant pour application $t$ la multiplication par $z_2$, on obtient une déformation d'ordre $k$ de l'espace vectoriel échelonné $\Et/t\Et$.
 En effet,
 $$N^0_s( t)=\sup_{| f|_s \leq 1} | t f|_s\leq |z_2|_s =s^2  $$
 et on a l'égalité lorsque $f$ est identiquement égale à $1$ donc, dans ce cas, $N^0_s( t) = s^2$. 
%%%%%%%%%%%%%%%%%%%%%%%%%%%%%%%%%%%%%%%%%%%%%%%
\subsection{Le groupe des isomorphismes}  
 Soit $(\Et,t)$ une déformation d'un espace vectoriel échelonné $E$
$$\ 0 \to \Et \stackrel{t}{\to} \Et \to E \to 0 \quad (*) .$$
Un morphisme de $\Et$ est appelé un {\em isomorphisme} s'il induit l'identité sur $E$ et s'il est inversible.
L'ensemble des isomorphismes forment le {\em groupe des isomorphismes} noté $\Isom(\Et)$.
 D'après le théorème \ref{T::produits}, l'exponentielle donne lieu à une application
 $$\Bt^1(E) \to \Isom(\Et),\ u \mapsto e^{tu}. $$

 L'exactitude de la suite $(*)$ implique que tout élément de  $\Isom(\Et)$ s'écrit de manière unique sous la forme $\Id+tu$, où $u$ est une application linéaire et $\Id \in \Lt(\Et)$ désigne l'application identité. En effet, pour $ \p \in \Isom(\Et),\ x \in \Et$, la projection de $ \p(a+tx)-\p$ sur $E$ est nulle, l'exactitude de la suite $(*)$ montre alors que $ \p(a+tx)-a$ est de la forme $t u(x)$. %En général, l'application $u$ est elle-même un morphisme~: 
%\begin{proposition} Si pour tout $s \in ]0,S[$, la quantité $\inf_{| x |_s=1}| t(x)|_s$ est non nulle alors tout isomorphisme de $\Et$ est de la forme
%$\Id+tu$ avec $u \in \Lt(\Et)$.
%\end{proposition}
%\begin{proof}
% Pour tout $s' \in ]0,S[$, il existe $s \in ]0,S[$ tel que l'application $\p$ définisse, par restriction, une application linéaire continue de $\Et_{s'}$ dans $\Et_s$. %On  a
%$$| t u(x) |_{s'} \geq (\inf_{| x |_{s'}=1}| tx|_{s'}) | u(x) |_s. $$
%Le membre de gauche tend vers $0$ quand $x$ tend vers $0_{\Et}$ dans $\Et_{s'}$ car $tu$ définit une application linéaire continue de $\Et_{s'}$ dans $%%une application linéaire continue de $\Et_{s'}$ dans $\Et_s$. La proposition est démontrée.
%\end{proof}    
%%%%%   %%%%%%%%%%%%%%%%%%%%%%%%%%%%%%%%%%%%%%%%
%%%%%%%%%%%%%%%%%%%%%%%%%%%%%%%%
 \subsection{Les théorèmes de stabilité}
\label{SS::groupes}
Dans le cadre des déformations d'espace vectoriel échelonnés, les théorèmes de $M$-détermination et de transversalité admettent des formulations simplifiées.
 \begin{theoreme}
 Soit $\Et$ un espace vectoriel échelonné, $a \in \Et$, $\alg$ un sous-espace vectoriel de $\Bt^1(\Et)$, $G$ un sous-groupe fermé de $\Isom(\Et)$ contenant $\exp(\alg)$ et agissant sur $a+t\Et$. Si
 \begin{enumerate}[{ \rm i)}]
 \item l'action de $\alg$ est échelonnée en $a$ ;
 \item l'application  $\rho:\alg \to \Et,\ u \mapsto u \cdot a$ possède un inverse borné ;
 \end{enumerate}
 alors l'orbite de $a$ sous l'action de $G$ est égale à $a+t\Et$.
 \end{theoreme}

\begin{theoreme}
  Soit $\alg$ un sous-espace vectoriel de $\Bt^1(\Et)$, $a \in \Et$, $G$ un sous-groupe fermé de $\Isom(\Et)$ contenant $\exp(\alg)$ et agissant sur $a+t\Et$.  Supposons que l'action de $\alg$ soit $\Ft$-échelonnée en $a$  et qu'il existe pour un certain $k \geq 0$ une application  bornée 
   $$j:\Ft \mapsto \Bt^k(\Et/\Ft,\alg)$$ 
 telle que  $j(\a)$ soit un inverse de  
 $$ \alg \to \Et/\Ft,  u \mapsto u \cdot (a+t\a) $$
 alors $t\Ft$ est une transversale à l'orbite de $a$ sous l'action de $G$ dans $a+t\Et$.
 \end{theoreme}

Dans un cadre abstrait, les théorèmes de stabilité sont donc des cas particuliers des résultats de détermination finie.
 %%%%%%%%%%%%%%%%%%%%%%%%%%%%%%%%%%%%%%%%%%%%%
\section{Structures échelonnées  en géométrie analytique}
%%%%%%%%%%%%%%%%%%%%%%%%
\subsection{Rappels : fibres d'un faisceau}
 Soit $\Ft$ un faisceau en espaces vectoriels topologiques définit sur un espace topologique $X$.  On appelle {\em fibre du faisceau $\Ft$ en un compact $K \subset X$}, noté $\Ft_K$, l'espace vectoriel topologique
 $$\Ft_K=\underrightarrow{\lim}\, \G(U,\Ft) $$
 où $U$ parcourt l'ensemble des ouverts contenant $K$ ordonné par l'inclusion.
   
 Un élément de $\Ft_K$ est une section du faisceau $\Ft$ au voisinage de $K$, pour laquelle on oublie de préciser la taille du voisinage de $K$ sur laquelle
 elle est définie.  Prendre la limite directe revient donc à identifier deux sections qui sont égales sur un ouvert contenant $K$~:
 $$f \sim g \iff \exists U \supset K,\ f_{\mid U}=g_{\mid U}.  $$ 
 Dans le cas où $K$ est réduit à un point, on retrouve la notion de germe en un point. Nous parlerons donc de germes de fonctions holomorphes en un compact. C'est une notion classique (voir par exemple \cite{Cartan_BSMF}). 
 
 En munissant les espaces vectoriels $ \G(U,\Ft)$ de la topologie de la convergence compacte, on munit la limite directe $\Ft_K$ d'une topologie. 
 
 Lorsque $K$ admet une base de voisinage donnée par des compacts, l'espace topologique $\Ft_K$ peut s'obtenir comme limite directe d'espaces de Banach de la façon suivante.
Notons $\Vt$ l'ensemble des voisinages compacts de $K$, ordonné par l'inclusion.
Pour $K' \in \Vt$, on note $\Bt_X(K')$ l'espace vectoriel des fonctions continues sur $K'$ qui sont holomorphes dans l'intérieur de $K'$.  C'est une espace de Banach pour la norme~: 
$$\Bt_X(K') \to \RM,\ f \mapsto \sup_{z \in K'}| f(z)|.$$
L'espace vectoriel $\Ft_K$  est limite directe des $ \Bt_X(K')$~:
$$\Ft_K =\underrightarrow{\lim}\, \Bt_X(K'), \ K' \in \Vt. $$  

  %%%%%%%%%%%%%%%%
   \subsection{L'image d'une application linéaire de rang fini}
Rappelons qu'une application linéaire continue entre espaces vectoriels topologiques $u:E \to F$ est appelée  {\em stricte} si elle induit un isomorphisme
 d'espaces vectoriels topologiques entre $E/\Ker u$ et $\Im u$ \cite{Bourbaki_topologie,Bourbaki_EVT}.  Si $E$ est métrisable alors pour que $u$ soit stricte,  il suffit que $u$ soit d'image fermée (théorème de l'image ouverte).  
  
 Le but de ce n° est de démontrer la
\begin{proposition}
\label{P::strict}
Soit  $K \subset \CM^n$ un compact et $M,N$ deux  modules de type fini sur l'anneau $\Ot_{\CM^n,K}$. Toute application $\Ot_{\CM^n,K}$-linéaire de $M$ vers $N$ est stricte.  
\end{proposition}
 Il suffit de démontrer   la proposition pour $M=\Ot_{\CM^n,K}^p$ et $N=\Ot_{\CM^n,K}^q$.   

Commençons par démontrer le
   \begin{lemme}
   \label{L::algebre} Soit $A$ une algèbre topologique. Si pour tout $x \in A$ l'image de la multiplication par $x$ est stricte alors toute application $A$-linéaire $A^n \to A^p$ est également stricte.
 \end{lemme}
%Par ailleurs, il suffit de montrer que l'image de chaque composante est fermé, donc on peut supposer $q=1$. Toute forme $A$-linéaire s'écrit sous la forme
 %$$A^p \to A,\ (x_1,\dots,x_p) \mapsto \sum_{i=1}^p a_i x_i,\ a_i \in A. $$
\begin{proof} 
 Soit $u,v:E \to F$ deux morphismes stricts, je dis qu'alors
\begin{enumerate}[{\rm i)}]
\item $(u,v):E \to F \times F,\ x \mapsto (u(x),v(x))$ est strict ;
\item  $u+v:E \to F,\ x \mapsto u(x)+v(x)$ est strict.
\end{enumerate}
La restriction d'un morphisme strict à un sous-espace vectoriel fermé est à nouveau un morphisme strict.  Donc la
restriction à diagonale $\D$ de l'application
 $$w:E \times E \to F \times F,\ (x,y) \mapsto (u(x),v(y)) $$
 est stricte. Ce qui démontre i). 
 
 Considérons le morphisme strict
 $$s:F \times F \to F,\ (x,y) \mapsto x+y .$$
 La composée de deux morphismes stricts est stricte donc $s \circ w_{| \D}=u+v$ est strict. Ce qui démontre ii). 
 
 Comme toute application linéaire est une somme finie d'applications de rang 1, les affirmations i) et ii) entraînent immédiatement le lemme.
 \end{proof}

 Si $\Ft$ est un faisceau analytique cohérent sur un espace analytique $X$ et si $K$ est un compact de $X$ alors l'espace vectoriel $\Ft_K$ ne vérifie pas la propriété de Baire. Cet espace vectoriel n'est donc pas métrisable, à moins, bien sûr, que $X$ ne soit réduit à un point, auquel cas c'est un espace vectoriel de dimension finie.  Cependant, Köthe a démontré un théorème pour les espaces LF, qui implique que le théorème de l'image ouverte reste valable pour  ces espaces~\cite{Koethe_Banach}.  Il nous reste donc à montrer que  la multiplication par $a $ est
d'image fermée, pour tout $a \in \Ot_{\CM^n,K}$. 

Pour cela, considérons des suites de germes en $K$ de fonctions holomorphes   $(y_k), (x_k)$ avec $y_k=a x_k$. Il s'agit de prouver que si $(y_k)$  est convergente alors $(x_k)$ l'est également. Pour cela, il suffit de démontrer que, pour chaque droite complexe $L \subset \CM^n $, la restriction des $x_k$ à $L$ définit une suite convergente. On est ainsi ramené au cas $n=1$. Si la limite existe elle est unique, il suffit donc de vérifier la propriété localement.
 
 Si $K$ ne contient pas de zéro de $a$, alors la suite $(x_k)$ s'écrit sous la forme 
 $$x_k=\frac{y_k}{a},\ a(0) \neq 0$$
 et comme $(y_k)$ est convergente, $(x_k)$ l'est également.
 
Supposons que $a$ s'annule sur $K$. Comme les zéros d'une fonction holomorphe sont isolés et comme la propriété d'être holomorphe est locale, on peut  supposer que $K=\{ 0 \}$. Posons $a(z)=z^db(z)$ avec $b(0) \neq 0$, on a
 alors $y_k(z)=z^dw_k(z) $ avec $w_k(0) \neq 0$. Ce qui montre que la suite $(x_k)$ s'écrit sous la forme 
 $$x_k=\frac{w_k}{b},\ b(0) \neq 0.$$
 Je dis que la suite $w_k$ est convergente. La suite $(y_k)$ étant convergente, il existe un disque fermé de rayon $s$ centré en l'origine $D_s \subset \CM$,
 tel que 
 \begin{enumerate}[{\rm i)}]
\item  $(y_k) \subset B_{\CM}(D_s)$ ;
\item $b(z) \neq 0,\ \forall z \in D_s.$
\end{enumerate} 
La suite $(w_k)$ est de Cauchy dans $B_{\CM}(D_s)$. Cet espace   est complet, donc $(w_k)$ et par suite $(x_k)$ convergent dans $B_{\CM}(D_s)$.  Ceci achève la démonstration de la proposition. (En fait nous avons montré que les modules sur un anneau de type $\Ot_{X ,K}$ forment une catégorie abélienne).

Soit $I \subset  \Ot_{\CM^n,K}$ un idéal engendré par $f_1,\dots,f_k \in  \Ot_{\CM^n,K}$.
En appliquant la proposition à l'application
$$\Ot_{\CM^n,K}^k \to  \Ot_{\CM^n,K},\ (a_1,\dots,a_k) \mapsto \sum_{i=1}^k a_i f_i ,$$
on obtient le
\begin{corollaire} Tout idéal de l'espace vectoriel topologique $ \Ot_{\CM^n,K}$ définit un sous-espace vectoriel fermé.
\end{corollaire}
 %%%%%%%%%%%%%%%%%%%%%%%%%%%%%%%%%
\subsection{Recouvrements échelonnés}
\label{SS::echelon}
Soit $(K_s) \subset \CM^n,\ s \in [0,S]$ une famille de compacts. Munissons $\CM^n$ de coordonnées $(z_1,\dots,z_n)$ et
notons $P$ le polycylindre
$$P=\{ z \in \CM^n : | z_i | \leq 1,\ i=1,\dots,n \}. $$
\begin{definition}
La famille $(K_s)$ est appelée un recouvrement échelonné si c'est une famille croissante telle que pour tout point $z \in K_s$ le polydisque $z+\s P$ est contenu dans le compact  $K_{s+\s}$, pour tous $s \in [0,S[$ et $\s \in [0,S-s[$.
\end{definition}
Nous allons construire des recouvrements échelonnés de la façon suivante. Soit
 $$\psi:\CM^n \supset \Omega \to [0,S], S \in [0,+\infty[$$ une fonction propre de classe $C^1$  définie sur un ouvert $\Omega \subset \CM^n$
 dont la dérivée est bornée.  Munissons  $\CM^n$ de la norme $\max_{i =1,\dots,n}| \cdot |$ et soit $\| \cdot \|$ la norme d'opérateur dans $L(\CM^n,\CM)$. Soit  $M$ un majorant de la norme des dérivées de $\psi$~:
 $$\sup_{z \in \Omega} \| D\psi(z)\| \leq M  $$
\begin{lemme} {La famille de compacts $K=(K_s)$ avec $K_s:=\psi^{-1}([0,Ms])$ est un recouvrement échelonné.}
\end{lemme}
\begin{proof}
Soit $z \in K_s$, la formule de Taylor donne~:
$$\psi(z+\s \dt)=\psi(z)+(\int_{t=0}^{t=1} D\psi(z+t\s \dt) dt)\s\dt,\ \dt \in P .$$
On a bien~:
$$  \psi(z+\s\dt)  \leq   \psi(z+\s\dt)  +  \sup_{t \in [0,1]} \| D\psi(z +t\s \dt)\| \s \leq Ms+M\s,$$
ce qui démontre le lemme.
  \end{proof}
 Si un recouvrement $K=(K_s)$ peut-être définit comme dans le lemme, et si de plus la fonction $\psi$ est pluri-sousharmonique, nous dirons que $K$ est {\em un recouvrement de
 Stein.}  
 %%%%%%%%%%%%%%%%%%%%%%%
 \subsection{La structure échelonnée  $B(X,K)$}
\label{SS::Banach}

Soit $K=(K_s)$ un recouvrement échelonné dans $\CM^n$.
  L'espace vectoriel topologique $\Ot_{\CM^n,K_0}$ est alors échelonné par les espaces de Banach $  B_X(K_s)$~:
  $$\Ot_{\CM^n,K_0} = \underrightarrow{\lim}\,B_X(K_s).$$

 Supposons que $K$ soit de Stein.  Si  $X$ est un sous-espace analytique de $\CM^n$ contenant $K $,  le théorème $A$ de Cartan montre  que l'on peut choisir une présentation du module $\Ot_{X,K_0}$(voir \cite{Cartan,Cartan_BSMF})~:
$$\Ot_{\CM^n,K_0} \stackrel{C}{\to} \Ot_{\CM^n,K_0} \to  \Ot_{X,K_0} \to 0$$

La proposition \ref{P::strict} du § précédent montre que l'image de l'application $C$ est fermée. L'espace vectoriel $\Ot_{X,K_0}$ se voit ainsi muni  d'une structure d'espace vectoriel échelonné. Cet structure ne dépend que du choix du recouvrement $K$ et pas du choix de la présentation, car si
$$\Ot_{\CM^n,K_0} \stackrel{C'}{\to} \Ot_{\CM^n,K_0} \to  \Ot_{X,K_0} \to 0 $$
désigne une autre présentation alors $\Im C \cap B_{\CM^n}(K_s)=\Im C' \cap B_{\CM^n}(K_s)$. Nous noterons $B(X,K)$ l'espace vectoriel topologique
$ \Ot_{X,K_0}$ muni de cette structure échelonnée. 

(Afin d'éviter toute confusion avec l'étude des voisinages priviliégiés, on prendra soin de remarquer que le sous-espace vectoriel ${\rm Im}\, C \cap B_{\CM^n}(K_s)$ n'est pas toujours égal à l'image par $C$ de $B_{\CM^n}(K_s)$ \cite{Cartan_ideaux,Douady_these}.)

On définit de même les échelonnements des fibres en un compact pour un faisceau analytique cohérent $\Ft$   sur  $X$~: on choisit une présentation de
$\Ft_{K_0}$~:
$$\Ot_{\CM^n,K_0}^p \stackrel{C}{\to} \Ot_{\CM^n,K_0}^q \to  \Ft_{K_0} \to 0,$$
et on prend sur $\Ft$ la structure échelonnée induite de $B(X,K)$.

Un résultat dû à Grothendieck montre que toute suite de Cauchy dans $\Ot_{X ,K_0}$ est en fait contenue dans un des  $\Bt_X(K_s) $ (voir \cite{Grothendieck_EVT}, Chapitre 3, Partie 1, Théorème 1 ou Partie 3, Proposition 5).  Par conséquent, en géométrie analytique, les ensembles fermés que
au sens des espaces vectoriels échelonnés (Chapitre I, \ref{SS::morphismes}) coïncident avec ceux de la topologie forte.
  
 %%%%%%%%%%%%%%%%%

%%%%%%%%%%%%%%%%%%%%%%%%%%%%%%%
 \subsection{La structure échelonnée  $H(X,K)$}
Soit $X,K$ comme précédemment. Notons $U_s$ l'intérieur de $K_s$.
Les espaces de Hilbert
$$H_{\CM^n}(U_s):=\bigcup_s L^2(U_s,\CM) \cap \G(U_s,\Ot_{\CM^n}).$$
munissent  l'espace vectoriel $\Ot_{\CM^n,K_0}$ d'un échelonnement
$$\Ot_{\CM^n,K_0} = \underrightarrow{\lim}\,H_{\CM^n}(U_s), $$
que nous noterons $H(X,K)$.

Comme toute fonction continue sur un compact est intégrable, l'identité donne un morphisme $0$-borné
$$J :B(X,K) \to H(X,K) $$ 
\begin{proposition} L'inverse de l'application $J$ est un morphisme $1$-borné.
\end{proposition}
Pour $z \in U_s$ et $\s$ fixés, on pose
$$f(z+\s \dt)=\sum_{j \geq 0} a_j \s^j,\ a_j \in \CM^n. $$
Comme le recouvrement $K$ est échelonnée, on a
$$| f |_{s+\s}^2 \geq \int_{z+\s P} | f(z)|^2 dV=\sum_{j \geq 0} |a_j|^2 \s^{2j+2}  $$
car $z+\s P \subset K_{s+\s}$. On en déduit l'inégalité~:
$$ |f(z) | \leq  \s^{-1} | f |_{s+\s}.$$
 Ce qui démontre la proposition.  
%%%%%%%%%%%%%%%
%%%%%%%%%%%%%%%%%%%%%%%%
\subsection{Inégalités de Cauchy}
\label{SS::lemme}
Soit $K$ un recouvrement échelonné de Stein dans une variété analytique $X$.
 \begin{proposition}
 \label{P::Cauchy}  Tout opérateur linéaire aux dérivées partielles d'ordre $k$ définit une application
 $k$-bornée de $ B(X,K)$.
 \end{proposition}  
D'une part, la propriété d'être borné est conservée par passage au quotient, et d'autre part, tout opérateur linéaire aux dérivées partielles d'ordre $k$ est obtenu par composition et addition d'opérateurs d'ordre~$1$. Il suffit donc de démontrer la proposition pour $X=\CM^n$, $k=1$. 

Dans ce cas, je dis qu'on a les inégalités
 $$ |\d_{z_i}f|_s \leq\frac{1}{\s}| f |_{s+\s},\ i=1,\dots,n   $$
 Reprenons les notations du  \ref{SS::Banach}.  
  Soit $f \in B_{\CM^n}(K_{s+\s})$ et  $z \in K_s$. Après une intégration par partie, la formule de Cauchy donne
$$\d_{z_i}f(z)=\int_{\d(z+\s P)}\frac{f(z)}{(\xi_i-z_i)\prod_{i=1}^n(\xi_i-z_i)}d\xi_1 \dots d\xi_n,\ z=(z_1,\dots,z_n)$$
En prenant des coordonnées polaires $\xi_i=z_i+\s e^{\sqrt{-1} \theta_i}$, on obtient~:
$$ |\d_{z_i}f(z)| \leq \frac{1}{\s} \sup_{w \in z+\s P}| f(w)| .   $$
Comme $K$ est échelonné, on a l'inclusion~:
$$z+\s P \subset K_{s+\s},\ z \in K_s $$
donc
$$ \frac{1}{\s} \sup_{w \in z+\s P}| f(w)| \leq \frac{1}{\s}| f |_{s+\s}.$$
Ceci achève la démonstration de la proposition.
%%%%%%%%%%%
 %%%%%%%%%%%%%%%%
 \subsection{Cas réel}
Supposons la variété $X $ munie d'une involution anti-holomorphe 
$$\tau:X \to X $$
et soit $X_{\RM}$ le lieu des points fixes de cette involution. Nous dirons qu'une fonction analytique est réelle si
$$f(\tau z)=\overline{ f(z)},\ \forall z.$$
En particulier, le domaine d'holomorphie de $f$ doit être invariant par $\tau$. On note $\Rt_{X}$ le faisceau des fonctions analytiques réelles
sur $X$.

 Soit $K$ un compact contenu dans $X_{\RM }$.  L'algèbre $\Rt_{X,K }$ est un sous-espace vectoriel topologique fermé
 de~$\Ot_{X ,K}$. Elle hérite par conséquent des structure échelonnées de $\Ot_{X ,K}$.

La partie réelle d'une sous-variété $X \subset \CM^n$ de Stein définie par des fonctions analytiques
 $$g_1,\dots,g_n:X \to \CM ,$$ admet des recouvrements échelonnés de Stein.
Il suffit, en effet, de poser\footnote{Cette fonction m'a été suggérée par P. Dingoyan.}~: 
 $$ \psi=\sum_{i=1}^n | g_i|+\sum_{i=1}^n | z-\tau (z) |$$ 
et de considérer un recouvrement associé à $\psi$ comme dans \ref{SS::echelon}.  Ces recouvrements $K$ sont invariants par l'involution $\tau$.
 %%%%%%%%%%%%%%%%%%%%%%%%%%%%%%%%%%%%%%%%%%%%
%%%%%%%%%%%
 %%%%%%%%%%%%%%%%%%%
  \subsection{Produits de Hadamard}
 \label{SS::Hadamard}
Soit $A$ une algèbre sur un corps $k$ et $M$ un $A$-module libre. Supposons que $M$ et $A$ soient munis
de topologies compatibles avec les structures d'algèbre et de module. Soit $(e_i), i \in I$ une famille libre de $M$
telle que l'adhérence du module engendré par les $e_i$ soit $M$ tout entier.

\begin{definition}
Le produit de Hadamard de deux éléments $f:=\sum_{i \in I} a_i e_i$, $g:=\sum_{i \in I }b_i e_i \in M$, par rapport à la base $e_i$, est défini par la série
$$f \star g:= \sum_{i \in I} a_i b_i\, e_i.$$
 \end{definition}
   
Considérons à présent le cas $A=k=\CM$ et $I=\ZM^n$. Soit une variété $X \subset \CM^n$ et  $K=(K_s)$ un recouvrement échelonné dans $X$. Nous dirons que les $(e_i)$ forment une base orthogonale de $H(X,K)$  si les $e_i$ sont orthogonaux l'espace de Hilbert $H(X,K)_s$
et s'ils engendrent un sous-espace vectoriel dense,  pour tout $s$.

 \begin{proposition}
 \label{P::Diophante}
Soit $(e_i),\ i \in \ZM^n$, une base orthogonale de $H(X,K)$ et $\l \in \CM^n$ un vecteur $(C,\tau)$-diophantien.
Supposons qu'il existe une application $k$-bornée telle que $Le_i=\a_i e_i$ avec $ |\s(i)^\tau| \leq \a_i$ alors
  le produit de Hadamard par la fonction 
  $$g:z \mapsto \sum_{i \in \ZM^n \setminus \{ 0 \}} \frac{1}{(\l,i)} e_i$$ est une application
  $k$-bornée.
\end{proposition}
\begin{proof}
Comme la base $(e_i)$ est orthogonale, on a:
$$| f \star g |^2_s \leq \sum_{i\in \ZM^n \setminus \{ 0 \}} |a_i b_i| | e_i|_s^2+ |a_0 b_0| | e_0|_s^2,$$
$f=\sum a_ie_i$, $g=\sum b_i e_i$. 
Comme $\l$ est diophantien, il vient~:
$$| f \star g |^2_s \leq C^2  \sum_{i \in \ZM^n \setminus \{ 0 \}}  ( \s(i)^\tau |a_i|)^2|e_i|_s^2 + |a_0 b_0| | e_0|_s^2 \leq C^2\, | Lf |_s^2+|b_0|^2 | f |_s^2  $$
\end{proof}
Comme l'identité induit un morphisme $1$-borné~: 
$$I:H(X,K) \to B(X,K),$$
la proposition reste valable pour les échelonnements de $B(X,K)$, à condition de remplacer $k$ par $k+1$ dans la conclusion.

Voici deux exemples qui nous serons utile par la suite. 
\begin{exemple} Supposons que $K$ soit un recouvrement échelonné de l'origine dans $\CM^n$ munit de coordonnées $z_1,\dots,z_n$. L'opérateur aux dérivées partielles
$$L:=\Id+(z_1\d_{z_1})^\tau+\dots+ (z_n \d_{z_n})^\tau$$
est d'ordre $\tau$ et $K$ est un recouvrement
échelonnée donc $L$ est un opérateur  $(\tau+1)$-borné de $H(X,K)$. Si le vecteur $\l=(\l_1,\dots,\l_n)$ est $(C,\tau)$-diophantien alors
 le produit de Hadamard par la fonction 
 $$g:z \mapsto \sum_{i \in \NM^n} \frac{1}{(\l,i)} e_i$$ définit également un opérateur $(\tau+1)$-borné .
\end{exemple}
\begin{exemple}Considérons le tore complexe 
$$(\CM^*)^n=\{ (z_1,\dots,z_n):\ \forall i, z_i \neq 0 \}.$$ 
Ce tore est muni de l'involution anti-holomorphe
$$(z_1,\dots,z_n) \mapsto (\frac{1}{\bar z_1},\dots,\frac{1}{\bar z_n}). $$
Sa partie réelle  est un tore réel $T$ de dimension $n$. Soit $K$ un recouvrement échelonné de $T$. Comme dans l'exemple précédent, la proposition implique que  le produit de Hadamard par la fonction 
$$g:z \mapsto \sum_{i \in \ZM^n \setminus \{ 0 \}} \frac{1}{(\l,i)} e_i$$
 est borné pourvu que le vecteur $\l=(\l_1,\dots,\l_n)$ soit diophantien. Cette opération est un inverse à droite de la  dérivation
  $$\Ot_{(\CM^*)^n,T}  \to  \Ot_{(\CM^*)^n,T},\ f \mapsto \sum_{j=1}^n\l_j z_j \d_{z_j} f.$$
 
Rappelons quelques points de terminologie. Le vecteur $\l=(\l_1,\dots,\l_n)$ définit un point de l'espace projectif $\PM^{n-1}$ appelé {\em fréquence} du champ de vecteurs $\sum_{j=1}^n\l_j z_j \d_{z_j}$.  Le flot définit par ce champ de vecteurs sur le tore est appelé quasi-périodique et lorsque la fréquence est diophantienne, on dit que le tore est {\em diophantien}.
\end{exemple}
On peut adapter, sans difficulté, les considérations de ce n° aux cas où la fonction $g$ de la proposition \ref{P::Diophante} dépend de paramètres. 
\section{Applications}
%%%%%%%%%%%%%%
%%%%%%%%%%%%%%%%%%%%%
  %%%%%%%%%%%%%%%%%%%%%%
 
%%%%%%%%
\subsection{Retour sur Poincaré et Siegel}
Revenons aux orbites de l'action adjointe du groupe des automorphismes de $\Ot_{\CM^n,0}$ sur l'algèbre de Lie $\alg$ des dérivations de cet anneau. L'action infinitésimale en $v \in \alg$ est donnée par
$$w \mapsto [v,w]. $$
Il s'agit de démontrer que sous les hypothèses de Siegel, cette application admet un inverse borné. Supposons que $v$ soit  de la forme
$$v=\sum_{i=1}^n \l_i z_i \d_{z_i}.$$
\'Ecrivons $w$ sous la forme
$$w=\sum_{j \in \NM^n, i=1 }^n a_{ij}z^j \d_{z_i},\ j=(j_1,\dots,j_n). $$
Un calcul direct montre que le crochet est donné par~:
$$[v,w]=\sum_{j \in \NM^n, i=1}^n ((j,\l)-\l_i) a_{ij}z^j \d_{z_i}. $$
L'inverse de l'action infinitésimale est donné par 
$$j:\sum_{j \in \NM^n, i=1}^n b_{ij}z^j \d_{z_i} \mapsto \sum_{j \in \NM^n, i=1 }^n ((j,\l)-\l_i)^{-1} b_{ij}z^j \d_{z_i} .$$

Fixons une structure échelonnée $H(\CM^n,K)$ où $K$ est un recouvrement échelonné de l'origine.
Lorsque $0$ n'est pas contenu dans l'enveloppe convexe des $\l_i$ et si ceux-ci sont linéairement indépendants sur $\QM$ alors les quantités
$((j,\l)-\l_i)^{-1} $ restent bornées~\cite{Poincare_these} (voir également \cite{Arnold_edo}). L'application $j$ est donc $0$-bornée pour $H(X,K)$. Les considérations du n° précédent montrent que  lorsque $\l$ est $(C,\tau)$-diophantien, c'est une application $(n\tau+1)$-bornée.
Par conséquent, les théorèmes de Poincaré et de Siegel sont des conséquences du théorème de $M$-détermination.  
%%%%%%%%%%%%%%%%%%%%%%%%%%%%%%%%%%%%%%%%%%
\subsection{Tores invariants}
\label{SS::KAM}
Les tores invariants diophantiens d'un système hamiltonien intégrable sont stables par perturbation de l'hamiltonien, pourvu que l'hamiltonien soit isochroniquement non-dégénéré. C'est essentiellement le contenu du théorème KAM.  Nous allons examiner une variante  de ce résultat.  

Le fibré cotangent au tore
$$(\CM^*)^n=\{ (z_1,\dots,z_n):\ \forall i, z_i \neq 0 \}.$$
est trivialisable. Nous notons 
$$(z_1,\dots,z_n,\xi_1,\dots,\xi_n) $$
des coordonnées pour lesquelles la forme symplectique est donnée par
$$\frac{1}{\sqrt{-1}}\left(\frac{dz_1}{z_1} \w d\xi_1+ \frac{dz_2}{z_2} \w d\xi_2+\dots+ \frac{dz_n}{z_n} \w d\xi_n \right).$$
Si on pose $z_i=e^{\sqrt{-1}\p_i}$, elle admet la forme familière~:
$$d\p_1 \w d\xi_1+ d\p_2 \w d\xi_2+\dots+ d\p_n \w d\xi_n $$
 Nous utiliserons abusivement la même notation $T$ pour le tore réel de $(\CM^*)^n$ et pour le produit de ce tore par $\{ 0 \} \subset \CM$.
Nous désignerons par $X$ le produit du fibré cotangent au tore complexe avec la droite $\CM=\{ t\}$ et par $I$ l'idéal engendré par les $\xi_i$ aussi bien  dans $\Ot_{X,T}$ que dans $\Ot_{(\CM^*)^n,T}$.

La forme symplectique induit une structure de Poisson $\CM\{ t\}$-linéaire sur l'algèbre $\Ot_{X,T}$.
 
 Munissons $\Ot_{X,T}$ d'une structure échelonnée de type $H(X,T)$ et notons $\Et$ l'union des espaces orthogonaux au $\CM$-espace vectoriel de dimension $n$ engendré par $\xi_1,\dots,\xi_n$.
L'image par  l'application
 $$\Ot_{X,T} \to \Ot_{X,T},\ f \mapsto \int_T f \frac{dz_1 \w dz_2 \w \cdots \w dz_n}{z_1 \dots z_n}$$
 d'un élément de $\Et$ est égale à une constante modulo $I^2$
\begin{proposition}Soit $f \in \Ot_{X,T}$ une fonction de la forme
$$f=\sum_{i=1}^n \l_i \xi_i\ (\mod I^2),\ \l_i \in \CM\{ t \} $$
Si la fréquence $( \l_1:\dots: \l_n) \in \PM^{n-1}$ est indépendante de $t$ et si elle est diophantienne alors l'espace vectoriel $t\Ft$ avec
$$ \Ft:=I^2 \oplus \CM\{ t \}$$ est une transversale  à l'action du groupe des isomorphismes de Poisson  dans $f+t\Et$. 
\end{proposition} 

Pour démontrer la proposition, considérons l'action infinitésimale 
$$\Ot_{X,T} \mapsto \Ot_{X,T},\ h \mapsto \{ h,f \} $$
et identifions $\Ot_{X,T}$ à la sous-algèbre des dérivations hamiltoniennes exactes par l'application
$$\Ot_{X,T} \to \Der_{\Ot_{X,T}}(\Ot_{X,T}),\ h \mapsto \{  h,\cdot \}. $$
On construit un inverse $j(\a)$ bornée de l'action infinitésimale
$$\rho(\a):\Ot_{X,T} \to \Et/I^2 \mapsto   g \mapsto \{ g,f+\a \},\ \a \in I^2$$
en cherchant d'abord un inverse modulo $I$ puis en ajoutant une correction afin  d'obtenir l'inverse modulo $I^2$. 
Pour cela, commençons par remarquer que la base $\xi_1,\dots,\xi_n$ permet de scinder la suite exacte d'algèbres
$$0 \to I/I^2 \to \Ot_{X,T}/I^2 \to \Ot_{X,T}/I \to 0 .$$
On obtient ainsi un isomorphisme d'algèbres
$$ \Ot_{X,T}/I^2 \approx \Ot_{X,T}/I \oplus I/I^2.   $$
   Soit, à présent, $g$ la fonction définie par~:
$$g =\sum_{i \in \ZM^n \setminus \{ 0 \}} \frac{1}{(\l,i)} z^i$$
Les résultats du \ref{SS::Hadamard}  montrent que le produit de Hadamard par $g$ est borné et définit un inverse de l'action infintésimale modulo $I$.   
On définit l'inverse $j(\a)$ de $\rho(\a)$ par la formule~:
$$ \Ot_{X,T}/I \oplus (\Et \cap I)/(\Et \cap I^2):(m,n) \mapsto (m \star g,n \star g-\{ m \star g ,\a \} \star g). $$
L'application $j$  est  bornée, les conditions du théorème de transversalité sont donc satisfaites. La proposition est démontrée.
 
On a une proposition analogue dans un cadre analytique réel. Avant d'appliquer la proposition à un cas concret, rappelons la procédure de moyennisation~: à chaque fonction $H$, on peut associer la fonction
$$H_0=\left( \frac{1}{\sqrt{-1}}\right)^n\int_{T} H  \frac{dz_1 \w dz_2 \w \cdots \w dz_n}{z_1 \dots z_n}$$
qui ne dépend que des variables $\xi_1,\dots,\xi_n$~:
$$H_0=\sum_{i=1}^n \l_i \xi_i +\sum_{i,j=1}^n a_{ij} \xi_i \xi_j(\mod I^3). $$
La condition de {\em non-dégénérescence isochronique} est~: $\det (\int_T a_{ij}) \neq 0$.  Si cette condition est satisfaite, on peut trouver un germe de courbe
$$\a:(\CM,0) \to (\CM^n,0),\ t \mapsto (\xi_1(t),\dots,\xi_n(t)) $$
lisse
au dessus duquel $H_0$ définit un flot quasi-périodique de fréquence constante égale à $\l$. \`A l'aide de cette courbe, on construit une famille à un paramètre de translations
$$\xi \mapsto \xi-\a(t) $$
qui ramène $H$ à un élément de la forme $f+t\Et$ modulo $I^2$, avec $f$ comme dans la proposition.
La proposition montre que $f+tR$ est conjugué à $f$ par une isomorphisme de Poisson.  Il existe donc une courbe lisse dans $\CM^n$ au-dessus de laquelle le flot de $H$ est conjugué à celui de $H_0$ au-dessus de $\a$. Lorsque $H$ et $R$ sont réels, la partie réelle de cette déformation définit une famille à un paramètre de tores invariants du système hamiltonien de fréquence $\l$. C'est une variante du théorème de Kolmogorov~\cite{Kolmogorov_KAM}.

Voyons à présent la situation où $H_0$ est dégénérée, ce qui nous permettra de mieux comprendre les notions mises en jeu dans la conjecture de Herman. L'application
$$F:(\CM^n,0) \to \PM^{n-1},\ \xi \mapsto  [\d_{\xi_1}H_0:\d_{\xi_2}H_0:\dots:\d_{\xi_n}H_0]$$
n'étant plus lisse, la fibre au-dessus de $F(0)$ est une sous-variété complexe singulière  de dimension au moins 1. Dans le cas réel,
la fonction $H$ admet une famille de tores invariants paramétrés par la partie réelle de cette variété. Dans le pire des cas, cette variété peut-être réduite à un point (par exemple si son idéal contient la fonction $\xi_1^2+\xi_2^2+\dots+\xi_n^2$). 
 
 %%%%%%%%%%%%%%%%%%%%%%%%%%%%%%%%%%%%%%%%%%%%%%%%%
 \subsection{Version singulière de KAM}
 Considérons le système hamiltonien définit par le germe de fonction
 $$H:(\CM^{2n},0) \to (\CM,0),\ (q,p) \mapsto \sum_{i=1}^n \l_i q_i p_i $$
 dans $\CM^{2n}$ muni de la forme symplectique standard. Les variétés stables et instables sont deux plans 
 lagrangiens. Notons $I$ l'idéal engendré par les $q_ip_i$, $i=1,\dots,n$. Le théorème de transversalité donne lieu au résultat suivant~(l'application du théorème général à ce cas particulier est quasiment identique à celle de la proposition du n° précédent).
 \begin{proposition} Si le vecteur $\l=(\l_1,\dots,\l_n)$ est diophantien alors pour toute fonction de la forme
 $$H+R,\ R \in \Mt^3_{\CM^{2n},0} $$
 il existe un automorphisme symplectique $\p \in \Aut(\Ot_{\CM^{2n},0})$ tel que
 $$\p(H+R)=H \ (\mod I^2). $$
 \end{proposition}
 En particulier le flot de $H$ est linéarisable par un symplectomorphisme, sur une variété lagrangienne complexe
 symplectomorphe à la variété d'idéal $I$.
 %%%%%%%%%%%%%%%%%%%%%%%%%%%%%%%%%%%%%%%
 \subsection{L'état fondamental au voisinage d'un point critique}
 Passons à présent à une version quantique du résultat précédent. Pour cela, rappelons que l'algèbre de Heisenberg $\widehat \Ht$ est  la $\CM[[\hbar]]$-algèbre des séries formelles sur $2n$-générateurs vérifiant les relations
 $$[q_i,p_j]=\hbar \dt_{ij}. $$
 Tout élément de $f \in \widehat \Ht$ s'écrit sous la forme
 $$f=\sum_{i,j,k \geq 0}  a_{ijk}q^ip^j\hbar^k,\ q^ip^j:= q_1^{i_1}q_2^{i_1} \dots q_n^{i_n}p_1^{i_1} p_2^{i_2}\dots p_n^{i_n} $$
 avec $i=(i_1,\dots,i_n),\ j=(j_1,\dots,j_n)$. 
 On appelle transformée de Borel de $f$, la série formelle dans les variables commutatives $x,y$
 $$Bf:=\sum_{i,j,k \geq 0} a_{ijk}x^iy^j\frac{\hbar^k}{k!} $$
 Nous noterons $\Ht$ la sous-algèbre formée des séries dont la transformée de Borel est analytique~\cite{Pham_resurgence} (voir également \cite{quantique}). On note $\Mt^3_{\Ht}$ l'image réciproque de l'idéal maximal par $B$.  
  
 On pose $\Ht\{ t \}:=\Ht \hat \otimes_{\CM} \CM\{ t \}$ où $\hat \otimes$ désigne le produit tensoriel topologique (il n'y pas lieu de préciser lequel car $\CM\{ t \}$ est nucléaire) \cite{Grothendieck_PTT}. 
  \begin{proposition} Si le vecteur $\l=(\l_1,\dots,\l_n)$ est diophantien alors pour tout élément de la forme
 $$H+tR,  R \in \Mt^3_{\Ht} \hat \otimes_{\CM} \CM\{ t \}$$
 il existe un automorphisme de $\Ht\{ t \}$ tel que
 $$\p(H+tR)=H \ (\mod I^2). $$
 \end{proposition}
 En particulier la série perturbative du spectre de l'état fondamental est Borel-analytique (voir \cite{quantique} pour plus de détails).
%%%%%%%%%%%%%%%%%%%%%
     \bibliographystyle{amsplain}
\bibliography{master}

%%%%%%%%%%%%%%%%%%%%%%%%%%%%%%%%%
 \end{document}